\documentclass[tbtags,reqno]{amsart}

\usepackage{amsthm,amsopn,amsfonts,amssymb}

\edef\savecatcodeat{\the\catcode`@}
\catcode`\@=11

\def\tb@ifSpecChars#1#2{#1}
\def\tb@ifNoSpecChars#1#2{#2}

\def\tableau{%
  \bgroup
  \@ifstar{\let\Tif\tb@ifNoSpecChars\tb@tableauB}
          {\let\Tif\tb@ifSpecChars\tb@tableauB}}

\def\tb@tableauB{
  \@ifnextchar[{\tb@tableauC}{\tb@tableauC[]}}

\def\tb@tableauC[#1]{\hbox\bgroup%
    \let\\=\cr
    \def\bl{\global\let\tbcellF\tb@cellNF}%
    \def\tf{\global\let\tbcellF\tb@cellH}
    \dimen2=\ht\strutbox \advance\dimen2 by\dp\strutbox%
    \ifx\baselinestretch\undefined\relax%
    \else%
       \dimen0=100sp \dimen0=\baselinestretch\dimen0%
       \dimen2=100\dimen2 \divide\dimen2 by\dimen0%
    \fi%
    \let\tpos\tb@vcenter
    \tb@initYoung
    \tb@options#1\eoo
    \let\arrow\tb@arrow%
    \dimen0=\Tscale\dimen2%
    \dimen1=\dimen0 \advance\dimen1 by \tb@fframe%
    \lineskip=0pt\baselineskip=0pt
    \def\tb@nothing{}%
    \def\endcellno{$\rss\egroup\bss\egroup}
    \def\endcell{\endcellno\kern-\dimen0}
    \def\begincell{\vbox to\dimen0\bgroup\vss\hbox to\dimen0\bgroup\hss$}%
    \let\overlay\tb@overlay%
    \let\fl\tb@fl%
    \let\lss\hss\let\rss\hss\let\tss\vss\let\bss\vss
    \def\mkcell##1{
        \let\tbcellF\tb@cellD
        \def\tb@cellarg{##1}
        \ifx\tb@cellarg\tb@nothing\let\tb@cellarg\tb@cellE\fi%
        \begincell\tb@cellarg\endcellno
        \tbcellF}
    \let\savecellF\tbcellF
     \Tif{\catcode`,=4\catcode`|=\active}{}\tb@tableauD}%

\let\tb@savetableauD\tableauD
{
    \catcode`|=\active \catcode`*=\active \catcode`~=\active%
    \catcode`@=\active
\gdef\tableauD#1{%
  \Tif{
    \mathcode`|="8000 \mathcode`*="8000%
    \mathcode`~="8000 \mathcode`@="8000%
    \def@{\bullet}%
    \let|\cr
    \let*\tf
    \let~\sk
  }{}%
  \tpos{\tabskip=0pt\halign{&\mkcell{##}\cr#1\crcr}}%
  \global\let\tbcellF\savecellF
  \egroup
  \egroup}
}
\let\tb@tableauD\tableauD
\let\tableauD\tb@savetableauD
\let\tb@savetableauD\undefined

\def\tb@options#1{\ifx#1\eoo\relax\else\tb@option#1\expandafter\tb@options\fi}

\def\tb@option#1{%
  \if#1t\let\tpos\tb@vtop\fi
  \if#1c\let\tpos\tb@vcenter\fi
  \if#1b\let\tpos\vbox\fi
  \if#1F\tb@initFerrers\fi
  \if#1Y\tb@initYoung\fi
  \if#1s\tb@initSmall\fi
  \if#1m\tb@initMedium\fi
  \if#1l\tb@initLarge\fi
  \if#1p\tb@initPartition\fi
  \if#1a\tb@initArrow\fi
}

\def\tb@vcenter#1{\ifmmode\vcenter{#1}\else$\vcenter{#1}$\fi}

\def\tb@vtop#1{\hbox{\raise\ht\strutbox\hbox{\lower\dimen0\vtop{#1}}}}

\def\tb@initPartition{\def\Tscale{.3}}
\def\tb@initSmall{\def\Tscale{1}}
\def\tb@initMedium{\def\Tscale{2}}
\def\tb@initLarge{\def\Tscale{3}}

\def\tb@initArrow{\dimen2=1.25em}

\def\tb@initYoung{%
  \def\tb@cellE{}
  \let\tb@cellD\tb@cellN
  \def\sk{\global\let\tbcellF\tb@cellNF}}
\def\tb@initFerrers{%
  \def\tb@cellE{\bullet}
  \let\tb@cellD\tb@cellNF
  \def\sk{\bullet}}

\tb@initMedium

\def\tb@sframe#1{%
  \vbox to0pt{
    \vss
    \hbox to0pt{%
      \hss
      \vbox to\dimen1{
        \hrule depth #1 height0pt
        \vss
        \hbox to\dimen1{
          \vrule width #1 height\dimen1
          \hss
          \vrule width #1
          }%
        \vss
        \hrule height #1 depth 0in
        }%
      \kern-\tb@hframe
      }%
    \kern-\tb@hframe}}

\def\tb@hframe{.2pt}\def\tb@fframe{.4pt}\def\tb@bframe{2pt}
\def\tb@cellH{\tb@sframe{\tb@bframe}}       
\def\tb@cellNF{}                            
\def\tb@cellN{\tb@sframe{\tb@fframe}}       
\let\tbcellF\tb@cellN                       

\def\tb@rpad{1pt}
\def\tb@lpad{1pt}
\def\tb@tpad{1.8pt}
\def\tb@bpad{1.8pt}

\def\tb@overlay{\endcell\@ifnextchar[{\tb@overlaya}{\begincell}}
\def\tb@overlaya[#1]{\vbox to\dimen0\bgroup%
  \tb@overlayoptions#1\eoo%
  \tss\hbox to\dimen0\bgroup\lss$}
\def\tb@overlayoptions#1{\ifx#1\eoo\relax\else\tb@overlayoption#1\expandafter\tb@overlayoptions\fi}

\def\tb@overlayoption#1{
  \if#1t\def\tss{\vskip\tb@tpad}\let\bss\vss\fi
  \if#1c\let\tss\vss\let\bss\vss\fi
  \if#1b\def\bss{\vskip\tb@bpad}\let\tss\vss\fi
  \if#1l\def\lss{\hskip\tb@lpad}\let\rss\hss\fi
  \if#1m\let\lss\hss\let\rss\hss\fi
  \if#1r\def\rss{\hskip\tb@rpad}\let\lss\hss\fi
}

\def\tb@fl{\endcell\begincell\vrule depth 0pt width \dimen0 height \dimen0 \endcell\begincell}

\@ifundefined{diagram}{}{
\def\tb@arrowpad{.5}

\newoptcommand{\tb@arrow}{\@ne}[2]{%
  \endcell
   \begingroup%
   \let\dg@getnodesize\tb@getnodesize
   \dg@USERSIZE=#1\relax%
   \ifnum\dg@USERSIZE<\@ne \dg@USERSIZE=\@ne \fi%
   \dg@parse{#2}%
   \dg@label{\tb@draw{#1}{#2}}}

\def\tb@getnodesize#1#2#3#4#5{\dimen3=\tb@arrowpad\dimen2 #4=\dimen3 #5=\dimen3\relax}
\def\tb@getnodesize#1#2#3#4#5{\ifnum#2=0\ifnum#3=0\tb@getnodesizetail{#4}{#5}\else\tb@getnodesizehead{#4}{#5}\fi\else\tb@getnodesizehead{#4}{#5}\fi}
\def\tb@getnodesizetail#1#2{\dimen3=.5\dimen2 #1=\dimen3 #2=\dimen3}
\def\tb@getnodesizehead#1#2{\dimen3=.5\dimen2 #1=\dimen3 #2=\dimen3}

\def\tb@draw#1#2#3#4{%
        \dg@X=0\dg@Y=0\dg@XGRID=1\dg@YGRID=1\unitlength=.001\dimen0%
        \dg@LBLOFF=\dgLABELOFFSET \divide\dg@LBLOFF\unitlength%
        \dg@drawcalc
        \begincell
        \let\lams@arrow\tb@lams@arrow
        \begin{picture}(0,0)\begingroup\dg@draw{#1}{#2}{#3}{#4}\end{picture}%
        \endcell
        \endgroup
        \begincell}
}

\def\tb@lams@arrow#1#2{%
 \lams@firstx\z@\lams@firsty\z@
 \lams@lastx#1\relax\lams@lasty#2\relax
 \lams@center\z@
 \N@false\E@false\H@false\V@false
 \ifdim\lams@lastx>\z@\E@true\fi
 \ifdim\lams@lastx=\z@\V@true\fi
 \ifdim\lams@lasty>\z@\N@true\fi
 \ifdim\lams@lasty=\z@\H@true\fi
 \NESW@false
 \ifN@\ifE@\NESW@true\fi\else\ifE@\else\NESW@true\fi\fi
 \ifH@\else\ifV@\else
  \lams@slope
  \ifnum\lams@tani>\lams@tanii
   \lams@ht\ten@\p@\lams@wd\ten@\p@
   \multiply\lams@wd\lams@tanii\divide\lams@wd\lams@tani
  \else
   \lams@wd\ten@\p@\lams@ht\ten@\p@
   \divide\lams@ht\lams@tanii\multiply\lams@ht\lams@tani
  \fi
 \fi\fi
 \ifH@  \lams@harrow
 \else\ifV@ \lams@varrow
 \else \lams@darrow
 \fi\fi
}

\catcode`\@=\savecatcodeat
\let\savecatcodeat\undefined

\numberwithin{equation}{section}

\makeatletter
\renewcommand{\subsubsection}{\@startsection
{subsubsection}
{3}
{0mm}
{\baselineskip}
{-0.5\baselineskip}
{\normalfont\normalsize\bfseries}}
\makeatother

\newtheorem{theorem}{Theorem}
\newtheorem{lemma}[theorem]{Lemma}
\newtheorem{proposition}[theorem]{Proposition}

\newtheorem{corollary}[theorem]{Corollary} 

\newtheorem{definition}[theorem]{Definition}

\newtheorem{property}[theorem]{Property}

\theoremstyle{remark}

\newtheorem*{acknow}{Acknowledgments}

\def\sas{\smallskip}

\def\charge{ {\rm {charge}}}
\def\cocharge{ {\rm {cocharge}}}
\def\shape{ {\rm {shape}}}
\def\stat { {\rm {Stat}}}
\def\dv { {\rm {dv}}}
\def\tab  { \mathsf T }
\def\sstab  {  T }
\def\stab  { \mathcal T }
\def \Sym {\mathcal Sym}

\begin{document}

\title[Tableaux statistics for $2$ part Macdonald polynomials]
{Tableaux statistics for two part Macdonald polynomials }

\author{L. Lapointe \and  J. Morse }

\address
{University of California, San Diego, 9500 Gilman Drive, La Jolla, CA
92093, USA 
}
 
\begin{abstract}
The Macdonald polynomials expanded in terms of
a modified Schur function basis have coefficients
called the $q,t$-Kostka polynomials.
We define operators to build standard tableaux and show
that they are equivalent to creation operators that recursively
build the Macdonald polynomials indexed by two part partitions.
We uncover a new basis for these particular Macdonald polynomials
and in doing so are able to give an explicit description of 
their associated $q,t$-Kostka coefficients by 
assigning a statistic in $q$ and $t$ to each standard tableau.
\end{abstract}

\maketitle

\section{Introduction }

The Macdonald polynomials, $J_\lambda[X;q,t]$,
are a two parameter family of polynomials
in $N$ variables, forming a basis for the space of symmetric functions.
The polynomials,  expanded in a modified Schur function basis
$\{S_\lambda[X^t]\}_{\lambda}$,
have coefficients called the
$q,t$-Kostka polynomials, $K_{\lambda\mu}(q,t)$.
Macdonald has conjectured that $K_{\lambda \mu}(q,t)$ is
a polynomial in $q$ and $t$ with positive integer coefficients.  
The associated representation-theoretic model \cite{[GH]} and computer
experimentation strongly suggest the existence of a statistic 
in $q$ and $t$ on standard tableaux which would 
validate this conjecture.  The case where $\mu$ is a partition 
with no more than 2 parts has been shown \cite{[Fish]} 
using rigged configurations which are in bijection with standard 
tableaux.  We present a statistic for this special case, 
obtained directly from standard tableaux,
by defining operators on standard tableaux
that correspond to
creation operators that build the Macdonald 
polynomials recursively.

The paper is divided as follows:  section two covers basic 
definitions used in symmetric function theory.  
We present in the third section,
the creation operator $B_2\!=\!tB_2^{(0)}\!+\!B_2^{(1)}q^{-D-1}$,
showing several properties that include the action of $B_2^{(0)}$ on the 
Hall-Littlewood polynomials $H_{\lambda}[X;q,t]$ 
and an expansion of $J_{\lambda}[X;q,t]$ in terms of products 
of $B_2^{(0)}$ and $B_2^{(1)}$ for $\ell(\lambda) \leq 2$.  
Further, a new basis for these Macdonald polynomials 
having coefficients in the parameters $q$ and $t$ with positive 
integer coefficients is uncovered.
The fourth section begins with basic definitions 
used in tableaux theory and then introduces 
operators on tableaux which correspond,
under a morphism $\digamma$, to the operators of section 3.  
Finally, a statistic on standard tableaux is 
presented in the fifth section.

\smallskip

\section{Definitions }
Symmetric polynomials are indexed by partitions, or decreasing sequences of
non-negative integers $\lambda =(\lambda_1,\lambda_2,\ldots)$ such that
$\lambda_1 \ge \lambda_2 \ge \dots$. 
The order of $\lambda$ is $|\lambda| = \lambda_1 + \lambda_2 + \dots$, 
the number of non-zero parts in $\lambda$ is denoted $\ell(\lambda)$ and  $n(\lambda)$ stands for
$\sum_{i} (i-1) \lambda_i$.
We define a {dominance order} on partitions 
where, given two partitions such that $|\lambda|=|\mu|$, 
$\lambda\leq\mu$ when $\lambda_1+\cdots+\lambda_i\leq
\mu_1+\cdots+\mu_i$ for all $i$.  
A partition $\lambda$ may be associated to 
a Young diagram with
$\lambda_i$ lattice squares in the $i^{th}$
row, from the bottom to top.  For instance, the Young diagram of $\lambda=(5,4,2,1)$ is
\begin{equation}
\tiny{\tableau*[sbY]{ \cr & \cr & & & \cr & & & & }} \, .
\end{equation}
For each square $s=(i,j)$ in the diagram of $\lambda$, let
$\ell'(s), \ell(s), a(s)$ and $a'(s)$ be respectively the number of
squares in the diagram of $\lambda$ to the south, north, east and west
of the square $s$.
The transposition of a Young diagram associated to $\lambda$ 
with respect to the main diagonal gives the conjugate partition
$\lambda'$.  For example, the conjugate of (2.1) is
\begin{equation}
\tiny{\tableau*[sbY]{ \cr & \cr & \cr & & \cr & & &  }} \, ,
\end{equation}
which gives $\lambda'=(4,3,2,2,1)$.  
\smallskip

We shall use $\lambda$-rings, needing only the formal ring of symmetric
functions $\Sym$ to act on the ring of rational functions in $x_1,\dots,x_N,q,t$,
with coefficients in $\mathbb R$.  The ring $\Sym$ is generated by power sums $\Psi_i$, $i=1,2,3\dots$.
The action of $\Psi_i$ on a rational function $\sum_{\alpha} c_{\alpha} u_{\alpha}/ \sum_{\beta} d_{\beta} v_{\beta}$
is by definition
\begin{equation}
\Psi_{i} \left[ \frac{\sum_{\alpha} c_{\alpha} u_{\alpha} 
 }{ \sum_{\beta} d_{\beta} v_{\beta} } \right]
 =\frac{\sum_{\alpha} c_{\alpha} u_{\alpha}^i }{ \sum_{\beta} d_{\beta} v_{\beta}^i}, 
\end{equation}
with $c_{\alpha},d_{\beta} \in \mathbb R$ and $u_{\alpha}, v_{\beta}$ 
monomials in $x_1,\dots,x_N,q,t$.  Since any symmetric
function is uniquely expressed in terms of the 
power sums, (2.3) extends to an action of
$\Sym$ on rational functions.  In particular, a symmetric function 
$f(x_1,\dots,x_N)$  can be denoted $f[x_1+\cdots+x_N]$.  We shall use 
the elements $X := x_1+\cdots+x_N$,  $X^{tq} :=
X(t-1)/(q-1)$ and $X^t := X(t-1)$. 
\smallskip

The Macdonald polynomials can now be defined.  We associate the number
\begin{equation} \label{1}
z_\lambda
        = 1^{m_1} m_1 !  \, 2^{m_2} m_2! \dotsm
\end{equation}
to a parition $\lambda$ with $m_i(\lambda)$ parts equal to $i$. 
Define a scalar product $\langle \ , \ \rangle_{q,t}$ 
on $\Sym \otimes \mathbb Q[q,t]$ by
\begin{equation}
\langle \Psi_\lambda[X], \Psi_\mu[X] \rangle_{q,t}
        =\delta_{\lambda \mu } z_\lambda \prod_{i=1}^{\ell(\lambda)} \frac{1-q^{\lambda_i}}{1-t^{\lambda_i}},
\end{equation}
where $\ell(\lambda)$ is the number of parts of $\lambda$.  The
Macdonald polynomials $J_\lambda [X; q,t]$ are uniquely specified by \cite{[Ma]}
\begin{align}
  \mathrm{(i)} \ &  \langle J_\lambda, J_\mu \rangle_{q,t} = 0, \qquad \text{if } \lambda \ne \mu , \\
  \mathrm{(ii)} \ &  J_\lambda[X;q,t] = \sum_{\mu \le \lambda} v_{\lambda\mu }(q,t) 
S_\mu[X] , \\
  \mathrm{(iii)} \ & v_{\lambda\lambda}(q,t)= 
\prod_{s \in \lambda} (1-q^{a(s)} t^{\ell (s)+1}),
\end{align}
where $S_{\mu}[X]$ is the usual Schur function.
\smallskip

It has been conjectured \cite{[Ma]} that the Macdonald polynomials
$J_{\lambda}[X;q,t]$, when expanded in the $\{S_{\mu}[X^t] \}_{\mu}$ basis, 
have expansion coefficients that are
polynomials in $q,t$ with positive integer coefficients. That is
\begin{equation}
J_{\lambda}[X;q,t] = \sum_{\mu} K_{\mu \lambda} (q,t) S_{\mu}[X^t] \, ,
\end{equation} 
with $K_{\mu \lambda} (q,t) \in \mathbb N[q,t]$.  
This conjucture has been proven \cite{[Fish]} in the
case where $\ell(\lambda) \leq 2$ (in fact, in the equivalent case where
$\lambda_1 \leq 2$) using rigged configurations.  
In this paper we present a description of 
these $K_{\mu\lambda}(q,t)$ with $\ell(\lambda)\leq 2$ in terms of
tableaux combinatorics.

\section{ Algebraic Side}

A Macdonald polynomial indexed by any partition
can be constructed by repeated application of 
creation operators $B_k$ defined \cite{[KN],[LV]} such that
\begin{equation}
B_k\, J_{\lambda}[X;q,t]
=
J_{\lambda+1^k}[X;q,t] \, , \qquad  \qquad \ell(\lambda) \leq k \, ,
\end{equation}
or more specifically, when $k=2$,
\begin{equation}
B_2\, J_{m,n}[X;q,t]
=
J_{m+1,n+1}[X;q,t] 
\, .
\end{equation}
If $\mathcal V$ is the $\mathbb Q[q,t]$-linear span 
of $\bigl\{ J_\lambda[X;q,t] \bigr\}_{\ell(\lambda) \leq 2}$, since 
it is known \cite{[Ma]} that 
\begin{equation}
J_{\lambda}[X;q,t]=
\sum_{ \mu \geq \lambda}c_{\mu \lambda} (q,t)\,S_\mu[X^{tq}] \, ,
\end{equation}
for some $c_{\mu \lambda}(q,t)\in\mathbb Q[q,t]$, 
$\mathcal
V$ is also the $\mathbb Q[q,t]$-linear span of 
$\bigl\{ S_\lambda[X^{tq}] \bigr\}_{\ell(\lambda) \leq 2}$.  
The action of $B_2$ on $\mathcal V$ can thus 
be defined by its action on 
$\bigl\{ S_\lambda[X^{tq}] \bigr\}_{\ell(\lambda) \leq 2}$,
introduced in \cite{[LLM]} as 
\begin{equation}
B_2 \, S_{m,n}[X^{tq}]= 
\hbox{ det } 
\left|
\begin{matrix}
(1-tq^{m+1})S_{m+1}[X^{tq}] & (1-q^{m+2})S_{m+2}[X^{tq}] \cr
(1-tq^n)S_{n}[X^{tq}] & (1-q^{n+1})S_{n+1}[X^{tq}] 
\end{matrix}
\right|\, . 
\end{equation}

It will be convenient to split the operator $B_2$ into a 
sum of two operators, 
\begin{equation}
B_2 = tB_2^{(0)}+B_2^{(1)} q^{-D-1},
\end{equation}
where $D$ is the operator 
such that $Df[X]=d\,f[X]$ on any homogeneous function of degree
$d$ and where $B_2^{(0)}$ and $B_2^{(1)}$ are defined 
on $\mathcal V$ by 
\begin{equation}
B_2^{(0)}\,S_{m,n}[X^{tq}] = 
\hbox{ det } 
\left|
\begin{matrix} 
-q^{m+1}S_{m+1}[X^{tq}] & (1-q^{m+2})S_{m+2}[X^{tq}] \cr
-q^nS_{n}[X^{tq}] & (1-q^{n+1})S_{n+1}[X^{tq}] 
\end{matrix}
\right|\, , 
\end{equation}
and
\begin{equation}
B_2^{(1)}\,S_{m,n}[X^{tq}] = 
q^{m+n+1} \hbox{ det } 
\left|
\begin{matrix} 
S_{m+1}[X^{tq}] & (1-q^{m+2})S_{m+2}[X^{tq}] \cr
S_{n}[X^{tq}] & (1-q^{n+1})S_{n+1}[X^{tq}] 
\end{matrix} 
\right|\, . 
\end{equation}
These expressions, obtained by expanding (3.4), provide that 
(3.5) holds on $\mathcal V$.

We use the following deformation of the Hall-Littlewood polynomials 
\begin{equation}
H_{\lambda}[X;q,t] = 
\sum_{\mu} q^{n(\lambda')} K_{\mu \lambda} (1/q,0) S_{\mu'}[X^t] \, ,
\end{equation}
which is a basis for the ring of symmetric functions since 
the $H_{\lambda}[X;q,t]$'s are linearly independent; i.e.,
the $K_{\mu \lambda}(1/q,0)$ matrix is triangular with respect to the partial
ordering and $K_{\lambda \lambda}(1/q,0)=1$.
These polynomials are specializations of the Macdonald polynomials. 
More precisely,  

\begin{proposition}
$H_{\lambda}[X;q,t]$ is obtained by taking the coefficient of
the maximal $t$-power, $t^{n(\lambda)}$, in the Macdonald polynomial $J_\lambda[X;q,t]$ 
expanded in $\bigl\{ S_\lambda[X^{t}] \bigr\}$ or 
$\bigl\{ S_\lambda[X^{tq}] \bigr\}$.
\begin{equation}
\begin{split}
H_{\lambda}[X;q,t] & := J_{\lambda}[X;q,t] \Big|_{t^{n(\lambda)}}^{\{t\}} 
= \sum_{\mu } K_{\mu \lambda} (q,t) \Big|_{t^{n(\lambda)}} 
S_{\mu}[X^t]  \\
& := J_{\lambda}[X;q,t] \Big|_{t^{n(\lambda)}}^{\{tq\}} 
=\sum_{\mu \geq \lambda }c_{\mu \lambda} (q,t) \Big|_{t^{n(\lambda)}}  \,S_\mu[X^{tq}] \, ,
\end{split}
\end{equation}
where $c_{\mu \lambda}(q,t)$ is as in (3.3) and 
the superscript indicates in which basis $J_\lambda$ is expanded.  
\end{proposition}
\noindent {\it Proof.} \quad  
We prove the expansion in the $\{ S_{\mu}[X^t]\}$ basis and the
other case follows from the relation:
\begin{equation}
S_{\lambda}[X^{t}] = \sum_{\mu} \bar v_{\mu \lambda}(q) S_{\mu} [X^{tq}]\,,
\end{equation}
for some $\bar v_{\mu\lambda}(q) \in \mathbb Q[q]$.  
We have $K_{\mu \lambda}(q,t)=q^{n(\lambda')}t^{n(\lambda)}K_{\mu'\lambda}(1/q,1/t)$ \cite{[Ma]}, 
which gives
\begin{equation}
K_{\mu \lambda} (q,t) \Big|_{t^{n(\lambda)}} 
= q^{n(\lambda')} K_{\mu' \lambda} (1/q,1/t) \Bigl|_{t^0} \, .
\end{equation}
Since $K_{\mu' \lambda} (1/q,1/t)$ 
is a polynomial in $1/q,1/t$ (see for instance \cite{[KN]} or \cite{[LV]}), (3.11) can be
rewritten as 
\begin{equation}
K_{\mu\lambda}(q,t) 
\Big|_{t^{n(\lambda)}}= q^{n(\lambda')} K_{\mu' \lambda} (1/q,0) \, ,
\end{equation}
which shows that the first equality in (3.9) is equivalent to (3.8).
\hfill $\square$
\sas
\begin{property}
We have
\begin{equation}
J_m[X;q,t] = H_m[X;q,t] = (q;q)_m S_m[X^{tq}] \, ,
\end{equation}
where the $q$-shifted factorial $(q;q)_m$ is such that 
\begin{equation}
(q;q)_m = (1-q)(1-q^2) \cdots (1-q^m) \, ,\qquad m>0\, ; \qquad (q;q)_0=1 \, .
\end{equation}
\end{property}

\noindent{\it Proof:}  It is shown in \cite{[Ma]} that 
$J_m[X;q,t] = (q;q)_m S_m[X^{tq}]$. 
Thus, from (3.9), $H_m[X;q,t]=(q;q)_m\Big|_{t^0}S_m[X^{tq}]$ 
is also equal to $J_m[X;q,t]$. \hfill $\square$
\smallskip

The characterization of $H_{\lambda}$ as a specialization of 
$J_\lambda$ leads to the following lemma;
\begin{lemma}
$\mathcal V$ is the 
$\mathbb Q[q,t]$-linear span of 
$\bigl\{ H_\lambda[X;q,t] \bigr\}_{\ell(\lambda) \leq 2}$ \, .
\end{lemma}

\noindent{\it Proof.}
(3.9) gives that the $\mathbb Q[q,t]$-linear span of 
$\bigl\{ H_\lambda[X;q,t] \bigr\}_{\ell(\lambda) \leq 2}$ 
is included in the $\mathbb Q[q,t]$-linear span of 
$\bigl\{ S_\lambda[X^{tq}] \bigr\}_{\ell(\lambda) \leq 2}$ 
which is equal to $\mathcal V$.  
The proposition then follows because 
$\bigl\{ H_\lambda[X;q,t] \bigr\}_{\ell(\lambda) \leq 2}$
is a linearly independent set.  \hfill $\square$
\smallskip

\noindent  The lemma implies that
$B_2^{(0)}$ can be defined on $\mathcal V$ by its action on $H_{m,n}[X;q,t]$
rather than by its action on $S_{m,n}[X;q,t]$.
\begin{property}  
The action of $B_2^{(0)}$ on $H_{m,n}[X;q,t]$ is given by
\begin{equation}
B_2^{(0)} H_{m,n}[X;q,t] = H_{m+1,n+1} [X;q,t] \, .
\end{equation}
\end{property}

\noindent{\it Proof.} \quad 
Definitions (3.6) and (3.7) show that the action of  
$B_2^{(0)}$ and $B_2^{(1)}$ on $S_{m,n}[X^{tq}]$ 
gives coefficients involving only the parameter $q$ when expanded in terms
of the $\bigl\{ S_\lambda[X^{qt}] \bigr\}_{\ell(\lambda) \leq 2}$ basis. 
This implies that the successive action of 
$B_2^{(0)}$ and $B_2^{(1)}$ on
$J_m[X;q,t]=(q;q)_m S_m[X^{tq}]$ produces coefficients involving 
only $q$ when expanded in the 
$\bigl\{ S_\lambda[X^{tq}] \bigr\}_{\ell(\lambda) \leq 2}$ basis.
We thus have, since
\begin{equation}
J_{m+\ell,\ell}[X;q,t]=
\left(B_2\right)^{\ell}
J_m[X;q,t]=\left(t B_2^{(0)}+B_2^{(1)}q^{-D-1}\right)^{\ell} J_m[X;q,t]
\, ,
\end{equation}
that the coefficient of the maximal $t$-power in
$J_{m+\ell,\ell}[X;q,t]$ is given by
\begin{equation}
J_{m+\ell,\ell}[X;q,t] \Big|_{t^{\ell}}^{\{tq\}} =
\Bigl( B_2^{(0)} \Bigr)^{\ell} J_m [X;q,t].
\end{equation}
The left side of this expression is 
$H_{m+\ell,\ell}[X;q,t]$ by (3.9) and $J_m[X;q,t]=H_m[X;q,t]$, 
hence (3.17) implies that (3.15) holds. \hfill $\square$
\smallskip

Another  property relating to the action of $B_2^{(0)}$ and $B_2^{(1)}$ 
can now be stated.
\begin{property} 
Let $\epsilon\in\{0,1\}$.  
$B_2^{(0)}$ and $B_2^{(1)}$ are such that 
\begin{equation}
(B_2^{(0)}+B_2^{(1)})^{m} H_\epsilon[X;q,t] = H_{2m+\epsilon}[X;q,t] \, .
\end{equation}
\end{property}

\noindent {\it Proof.} \quad Definitions (3.6) and (3.7) give that
\begin{equation}
(B_2^{(0)}+B_2^{(1)}) S_k[X^{tq}] = (1-q^{k+1})(1-q^{k+2}) S_{k+2}[X^{tq}] \, ,
\end{equation}
which, using $H_k[X;q,t] = (q;q)_k S_k[X^{tq}]$, implies that
\begin{equation}
(B_2^{(0)}+B_2^{(1)}) H_k[X;q,t] = H_{k+2}[X;q,t] \, .
\end{equation}
We apply this identity $m$ times, starting with $k=\epsilon$,
to complete the proof.\hfill $\square$ 
\smallskip

\noindent  With relation (3.5), this property shows that any
Macdonald polynomial indexed by $\leq 2$ part partitions can be
built using $B_2^{(0)}$ and $B_2^{(1)}$ since
$H_m[X;q,t]=J_m[X;q,t]$.
\smallskip

We have further that $B_2^{(0)}$ and $B_2^{(1)}$ satisfy a $q$-commutation relation.
\begin{property} 
On the space $\mathcal V$, we have 
\begin{equation}
B_2^{(1)} B_2^{(0)} = q \, B_2^{(0)} B_2^{(1)} \,. 
\end{equation}
\end{property}
\noindent {\it Proof.} \quad 
The action defined in (3.5) and (3.6), and the Pieri rule 
\begin{equation}
S_{m}[X^{tq}] S_{n}[X^{tq}] =\sum_{\ell=0}^{n} S_{m+n-\ell,\ell} [X^{tq}]\, , 
\quad\hbox{\rm where}\quad  m \geq n \, ,
\end{equation}
yield
\begin{equation}
\begin{split}
B_2^{(0)} B_2^{(1)} S_{m,n}[X^{tq}] & = q^{n+m+2} (q^{m+1}-q^n) 
\sum_{r = 0}^n \sum_{\ell=0}^{r}
(q^{r} -q^{m+n+3-r}) S_{m+n+4-\ell,\ell} [X^{tq}] \\
& \,\quad +\,  q^{2m+2n+5-\ell} (q^{m+1}-q^n) \sum_{\ell = 0}^n 
(q^{\ell+1} -1) S_{m+n+3-\ell,\ell+1} [X^{tq}] \\
&\,\quad + \, q^{n+m+1} (1-q^{n+1}) (q^{n+1}-q^{m+2}) \sum_{\ell = 0}^{n+1} 
S_{m+n+4-\ell,\ell} [X^{tq}] \\
& \,\quad + \, q^{2m+n+3} (1-q^{n+1}) (q^{n+2}-1) S_{m+2,n+2} [X^{tq}] 
\end{split}
\end{equation}
and 
\begin{equation}
\begin{split}
B_2^{(1)} B_2^{(0)} S_{m,n}[X^{tq}] & = q^{n+m+4} (q^{n}-q^{m+1}) \sum_{r = 0}^n \sum_{\ell=0}^{r}
(q^{m+n+3-r}-q^{r}) S_{m+n+4-\ell,\ell} [X^{tq}] \\
&\, \quad +\,  q^{n+m+3} (q^{n}-q^{m+1}) \sum_{\ell = 0}^n 
(1-q^{\ell+1}) S_{m+n+3-\ell,\ell+1} [X^{tq}] \\
& \,\quad +\,  q^{2m+n+5} (q^{n+1}-1) \sum_{\ell = 0}^{n+1} (q^{m+2}-q^{n+1}) 
S_{m+n+4-\ell,\ell} [X^{tq}] \\
&\, \quad + \, q^{2m+n+4} (1-q^{n+1}) (q^{n+2}-1) S_{m+2,n+2} [X^{tq}] \, .
\end{split}
\end{equation}
First exchange the order of summation in the right side of 
the top line in (3.23) and (3.24), and then use
\begin{equation}
\sum_{r =\ell}^n (q^r-q^{m+n+3-r})= 
\frac{(1-q^{n-\ell+1})(q^{\ell}-q^{m+3})}{(1-q)}.
\end{equation}
Next, send $\ell \to \ell\!-\!1$ 
in the right side of the second terms in (3.23) and (3.24).  
By fixing $\ell$, we have that 
$q B_2^{(0)}B_2^{(1)}S_{m,n}[X^{tq}]=B_2^{(0)} B_2^{(1)} S_{m,n}[X^{tq}]$ if
\begin{equation}
\begin{split}
&  \left(  \frac{(1-q^{n-\ell+1})(q^{\ell}-q^{m+3})}{(1-q)} +q^{m+n+4-\ell}
(q^{\ell}-1) +(q^{n+1}-1) \right) \\
& \qquad \qquad \qquad =   \left(  q \frac{(1-q^{n-\ell+1})(q^{\ell}-q^{m+3})}{(1-q)} +
(q^{\ell}-1) +q^{m+3}(q^{n+1}-1) \right) \, .
\end{split}
\end{equation}
This algebraic relation can be checked straightforwardly.  \hfill $\square$
\smallskip

\begin{definition} Let $v = (v_1,\dots,v_k)$ with all $v_i \in \{ 0,1\}$.  For
$\epsilon\in\{0,1\}$ we define
\begin{equation}
U_v^{(\epsilon)} = B_2^{(v_1)} \cdots B_2^{(v_k)} \cdot H_\epsilon[X;q,t] \, .
\end{equation}
\end{definition}
\begin{proposition} 
For $\epsilon,v_i\in\{0,1\}$ we have
\begin{equation}
J_{2 m+ \ell+\epsilon,\ell}[X;q,t] = 
\sum_{v=(v_1,\dots,v_{m+\ell})} q^{(1-d)  |v|_{\ell}+2 n(v)_{\ell}} 
t^{\ell -|v|_{\ell} } U_v^{(\epsilon)},  
\end{equation}
where $d=2m +2\ell+ \epsilon$, 
$|v|_{\ell}= v_1+ \cdots + v_{\ell}$ and 
$n(v)_{\ell}= v_2 + 2 v_3 + \cdots + (\ell-1)v_{\ell}$.  
\end{proposition}
\noindent {\it Proof.} \quad
From property 5, we have that
\begin{equation}
J_{2m+ \epsilon}[X;q,t] = \sum_{v=(v_1,\dots,v_{m})}
U_v^{(\epsilon)} \, ,
\end{equation}
where $v_i\in\{0,1\}$, proving (3.28) for $\ell = 0$.  
Proceeding by induction, we assume that (3.28) holds for every $\ell$. 
We thus have, acting with $B_2$, that
\begin{equation}
\begin{split}
& (t B_2^{(0)}+B_{2}^{(1)} q^{-D-1} )J_{2m+\ell +\epsilon,\ell}  
= \sum_{v'=(0,v)} q^{(1-d)  |v'|_{\ell+1}+ 2( n(v')_{\ell+1}-|v'|_{\ell+1})} t^{\ell+1 -|v'|_{\ell+1} } U_{v'}^{(\epsilon)}\\
& \qquad \qquad  +  
\sum_{v"=(1,v)} q^{(1-d)( |v''|_{\ell+1}-1)+
2( n(v'')_{\ell+1}-|v''|_{\ell+1}+1) -d-1} t^{\ell+1 -|v''|_{\ell+1} }
 U_{v''}^{(\epsilon)} \,.
\end{split}
\end{equation}
Combining the two sums, we obtain
\begin{equation} 
B_2 J_{2m+\ell +\epsilon,\ell} = 
\sum_{\bar v= (\bar v_1,\dots,\bar v_{m+\ell+1})} 
q^{\bigl( 1-(d+2) \bigr)  
|\bar v|_{\ell+1}+ 2 n(\bar v)_{\ell+1}} t^{\ell+1 -
|\bar v|_{\ell+1} } U_{\bar v}^{(\epsilon)},
\end{equation}
which completes the induction argument
since $B_2 J_{2m+\ell +\epsilon,\ell}= J_{2 m+\ell+1 +\epsilon,\ell+1}$. \hfill $\square$

\noindent{\it Example:}  
\quad We have
\begin{equation}
\begin{split}
J_{4,2}[X;q,t]&  = t^2 U_{0,0,0}^{(0)} + t^2 U_{0,0,1}^{(0)} 
+ tq^{-3} U_{0,1,0}^{(0)} + tq^{-3} U_{0,1,1}^{(0)} \\
& \quad +
 tq^{-5} U_{1,0,0}^{(0)} + tq^{-5} U_{1,0,1}^{(0)} 
+ q^{-8} U_{1,1,0}^{(0)} + q^{-8} U_{1,1,1}^{(0)} \, .
\end{split}
\end{equation}

\begin{corollary}  The maximal $t$-power in (3.28) is
\begin{equation}
H_{2m+\ell+\epsilon, \ell}[X;q,t] = \sum_{\bar v} U_{\bar v}^{(\epsilon)} \, ,
\end{equation}
summing over all $\bar v=(0^{\ell},v)$ where $v=(v_1,\dots,v_m)$ such that
$v_i \in \{0,1 \}$. 
\end{corollary}

\section{Tableaux Side}

\subsection{Definition and background}

Let $\mathcal A^{*}$ be the free monoid 
generated by the alphabet $\mathcal A = \{1,2,3, \dots \}$
and $\mathbb Q [\mathcal A^{*}]$ be the free algebra of $\mathcal A$.  
The elements of $\mathcal A^{*}$ are called words.  
The degree of a word $w$ is denoted $|w|$ and its 
image in the ring of polynomials $\mathbb Z [\mathcal A]$ is called 
the evaluation, denoted $ev(w)$.  For example, 
$w = 1 3 1  3 3 2$ has degree 6 and evaluation (2,1,3). 
A word $w$ of degree $n$ is said to be standard 
iff $ev(w) = (1,1,\dots,1)$.

A tableau $\tab$ will be the  pair $(\lambda,w)$, 
where $\lambda$ is a partition and $w$ is a word, 
such that $|\lambda|=|w|$.  We say that $\lambda$ is the shape of $\tab$.
A Young diagram associated to $\lambda$ filled with
the letters of $w$ from left to right and
top to bottom is a planar representation of  $\tab$.
For example, $\tab= ((4,2,2,1), 114356234)$ has 
the planar disposition 
\begin{equation}
\tab = \tiny{\tableau*[scY]{1 \cr 1&4 \cr 3&5 \cr 6&2&3&4 }} \, .
\end{equation} 
A semi-standard tableau is a tableau such that the entries in every rows
are nondecreasing  and such that the entries in every columns are increasing. 
In that case, we do not need to specify $\lambda$ in the pair  $(\lambda,w)$, since it
can be extracted from $w$.  For instance, $\sstab = 67~445~11123$ has
the representation 
\begin{equation}
\sstab = \tiny{\tableau*[scY]{
6 & 7 \cr
4 & 4 & 5 \cr
1 & 1 &1 &2 & 3 }} \, .
\end{equation}
Notice that a semi-standard tableau is a tableau.
Finally, a standard tableau $\stab$ is a 
semi-standard tableau of 
evaluation $(1,1,\dots,1)$.  For example, $\stab=7~46~1235 $ or
\begin{equation}
\stab = \tiny{\tableau*[scY]{
7  \cr
 4 & 6 \cr
1 & 2 &3 &5 }}  
\end{equation}
is a standard tableau.  

Obtained in the same manner as the transpose of a Young diagram,
we define $\tab^t$ to be the transpose of a tableau $\tab$.
For example, with $\tab$ as given in (4.1), we have
\begin{equation}
\tab^t = \tiny{\tableau*[scY]{ 4 \cr 3 \cr 2 & 5 &4 \cr 6 & 3 & 1 &1  }} \, .
\end{equation}
Notice that the transpositon of a 
semi-standard tableau may not be 
a semi-standard tableau, 
but the transposition of a standard tableau is always 
a standard tableau.

Words can be associated to numbers called the charge and cocharge
where
\begin{equation}
\charge (w) = n( ev(w)_P ) -\cocharge (w) \, ,
\end{equation}
for $ev(w)_P$ the partition obtained by reordering $ev(w)$.
The cocharge of a standard word $w$ is defined by the following algorithm;
\begin{enumerate}
\item[{1.}] Label the letter 1 in $w$ by $c_1 =0$ 
\item[{2.}] If the letter $i+1$ appears at the left of the letter $i$ in $w$, 
then $c_{i+1}=c_i+1$. Otherwise $c_{i+1}=c_i$. 
\item[{3.}] $\cocharge (w)= c_1 + \cdots + c_n$.
\end{enumerate}
For instance, $\cocharge(413265)= 0+0+1+2+2+3=8$.
Recall that  semi-standard tableaux and standard tableaux 
are simply words, and as such have an associated charge and cocharge. 

Lascoux and Sch\"utzenberger defined \cite{[LS1]} an action 
of the symmetric group on $\mathbb Z [\mathcal A^{*}]$
that sends a word of evaluation
$(ev_1,\dots, ev_i,ev_{i+1},\dots)$ 
to a word of evaluation 
$(ev_1,\dots, ev_{i+1},ev_{i},\dots)$
under an elementary transposition, $\sigma_i$. 
Their action induces the usual 
action of the symmetric group on $\mathbb Z [\mathcal A]$. 
For our purposes, we define this action only on words 
such that $(ev_i,ev_{i+1})\in\{(1,2),(2,1)\}$:
\begin{equation}
a a b \overset{\sigma_a}{\leftrightarrow} abb ,  \quad  
a b a  
\overset{\sigma_a}{\leftrightarrow} b b a, 
\quad b aa \overset{\sigma_a}{\leftrightarrow} bab\, ,
\end{equation}
where $a$ and $b$ stand for $i$ and $i+1$ respectively.  Furthermore, this action sends a semi-standard tableau to a 
semi-standard tableau while preserving its shape.
For example, $\sigma_4 (215345)= 215344$ and $ \sigma_3 \, 
\tiny{\tableau*[scY]{ 4 \cr 3 &5&5 \cr
1&2&3}} = \tiny{\tableau*[scY]{ 4 \cr 3 &5&5 \cr
1&2&4}}$.

We will use several simple linear operations on words. 
$\tau_k$ is a translation of $a$ to $a+k$ for 
every letter $a$ in an alphabet.  
For instance, $\tau_2(231567)=453789$.
The restriction of a word $w$ to the alphabet $a,b,c,\dots$,
is denoted $w_{\{a,b,c,\dots\}}$.  i.e., 
$w_{\{3,4\}}=3343$  for
$w= 1 2 3 3 4 2 23$. 
If $w$ is such that $w_{\{a,b\}}=ab$, 
$r_{(ab \to cd)}$ sends $ab$ to $cd$.  
For example, $r_{(23 \to 46)} (1 21 543)= 141546$.
We further define an operator $R_a$  
to remove all letters $a$ in $w$, and
finally, $A_{n+1,n+1}$ is an operator
on a tableau $\tab$, that adds 
a horizontal 2-strip of the boxes $n+1$ in all
the possible ways to $\tab$.  For instance
\begin{equation}
A_{5,5} \,  \,  \tiny{\tableau*[scY]{ 3 \cr 2 \cr 1&4}} = 
\tiny{\tableau*[scY]{ 5 \cr 3 \cr 2& 5 \cr 1&4}} + 
\tiny{\tableau*[scY]{ 5 \cr 3 \cr 2 \cr 1&4&5}} +
\tiny{\tableau*[scY]{ 3 \cr 2 &5 \cr 1&4&5}} +
\tiny{\tableau*[scY]{ 3 \cr 2 \cr 1&4&5&5}} \, .
\end{equation}

There exists a shape and cocharge preserving
standardization \cite{[LS1]} of semi-standard tableaux,
denoted by $VS$, that we will use only
on semi-standard tableaux with evaluation $(ev_1,\dots,ev_k)$, 
where $ev_i \in \{1,2 \}$.  In the case of a semi-standard tableau $\sstab$,
it is defined as
\begin{enumerate}
\item[1.] If $ev_1=2$ then $\sstab\to \tau_{1 } r_{(11 \to 01)} \sstab$.
Proceed to step 2.
\item[2.] If $ev(\sstab)=(1,1,\dots,1)$ then the standardization is complete.  
          Otherwise, 
          $\sstab\to\sigma_1 \cdots \sigma_{i-1} \sstab$
for the smallest $i$ such that $ev_i=2$. Proceed to step 1. 
\end{enumerate}
For example, 
$\sstab=45~235~124$ undergoes the following 
standardization process:
\begin{equation}
\tiny{\tableau*[scY]{
4 & 5 \cr
2 & 3 & 5  \cr
1 & 2 &  4 }} \, \to \,
\tiny{\tableau*[scY]{
4 & 5 \cr
2 & 3 & 5  \cr
1 & 1 & 4 }} \, \to \, 
\tiny{\tableau*[scY]{
5 & 6 \cr
3 & 4 & 6  \cr
1 & 2 & 5 }}
 \, \to \,
\tiny{\tableau*[scY]{
5 & 6 \cr
2 & 3 & 6  \cr
1 & 1 & 4 }}
\, \to \,
\tiny{\tableau*[scY]{
6 & 7 \cr
3 & 4 & 7  \cr
1 & 2 & 5 }}
 \, \to \, 
\tiny{\tableau*[scY]{
5 & 6 \cr
2 & 3 & 7  \cr
1 & 1 & 4 }}
\, \to \, 
\tiny{\tableau*[scY]{
6 & 7 \cr
3 & 4 & 8  \cr
1 & 2 & 5 }}
 \, .
\end{equation}
Note that $VS^{(n)}$ will denote the standardization of
only the first $n$ letters of a semi-standard tableau 
of degree $N\geq n$ which sends
such a tableau to a semi-standard tableau of evaluation 
$(1^n,ev_{n+1},\dots,ev_{N})$.

\subsection{Tableaux operators}

It is known \cite{[LS2]} that the Hall-Littlewood polynomials 
\begin{equation}
Q_{\lambda}[X;t] = \sum_{\mu} K_{\mu \lambda}(0,t) S_{\mu}[X^t],
\end{equation}
are equivalently expressed as a sum over semi-standard 
tableaux $\sstab$ such that 
\begin{equation}
Q_{\lambda}[X;t] = \sum_{\sstab ; ev(\sstab)=\lambda} 
t^{\charge(\sstab)} S_{\shape (\sstab )}[X^t].
\end{equation}
The substitution of $K_{\mu \lambda}(1/q,0)= K_{\mu' \lambda'}(0,1/q)$ \cite{[Ma]}
in expression (3.8) thus yields
\begin{equation}
\begin{split}
H_{\lambda}[X;q,t] & = \sum_{\mu} q^{n(\lambda')} K_{\mu' \lambda'}(0,1/q) S_{\mu'}[X^t] \\
  & = \sum_{\sstab ; ev(\sstab)=\lambda'} 
q^{\cocharge(\sstab)} S_{\shape (\sstab) }[X^t]\,,
\end{split}
\end{equation}
using $\cocharge(\sstab) = n(ev(\sstab))-\charge(\sstab)$.
Since standardization $VS$ preserves cocharge and shape, we then have
\begin{equation}
 H_{\lambda}[X;q,t] = \sum_{\sstab ; ev(\sstab)=\lambda'} 
q^{\cocharge(VS(\sstab))} 
S_{\shape \bigl(VS(\sstab) \bigr)}[X^t]\,.
\end{equation}
Further, defining a morphism $\digamma$ on semi-standard tableaux 
such that,
\begin{equation}
\digamma: \sstab \to q^{\cocharge(\sstab)} S_{\shape (\sstab ) }[X^t] \, ,
\end{equation}
we may express $H_\lambda[X;q,t]$ as the action of $\digamma$ on
a sum of semi-standard tableaux.
  More precisely, using the following definition:
\begin{definition}
For all partition $\lambda$ such that $\ell(\lambda) \leq 2$,
\begin{equation}
 \mathbb H_{\lambda} = \sum_{\sstab ; ev(\sstab)=\lambda'} VS(\sstab)
\, ,
\end{equation}
\end{definition}
\noindent we then have
\begin{equation}
\digamma(\mathbb H_{\lambda} )  =  H_{\lambda}[X;q,t]\, .
\end{equation}

We now define linear operators on standard tableaux. 
\begin{definition} 
On any standard tableau $\stab$ such that $|\stab|=n$, we define
\begin{equation}
 \mathbb B_2^{(0)}  : \stab \to VS( A_{n+1,n+1} \stab) = \tau_1 r_{(11 \to 01)} 
\sigma_1 \cdots \sigma_{n} A_{n+1,n+1} \stab
\end{equation} 
and
\begin{equation}
\mathbb B_2^{(1)} : \stab \to \Bigl( \mathbb B_2^{(0)} \stab^t \Bigr)^t \, .
\end{equation}
\end{definition}
\noindent Note that 
$\mathbb B_2^{(0)}$ and $\mathbb B_2^{(1)}$ 
send $\stab$ to a sum of standard tableaux.

\noindent {\it Example:} \quad 
Given a standard tableau $\stab=\tiny{{\tableau*[scY]{ 3 \cr 1&2&4 }}}\,$, 
we add all the possible horizontal 2-strips containing the letter $5$
and then standardize.
\begin{equation}
\begin{split}
\mathbb B_2^{(0)} 
\left( \, \tiny{{\tableau*[scY]{ 3 \cr 1&2&4 } }}\, \right) &  = 
VS \left( \, \tiny{{\tableau*[scY]{ 5 \cr 3 & 5 \cr 1&2&4 } }} 
+ \tiny{{\tableau*[scY]{ 5 \cr 3 \cr  1&2&4&5 } }}
+\tiny{{\tableau*[scY]{ 3 & 5 & 5 \cr 1&2&4 } }}
+\tiny{{\tableau*[scY]{ 3 & 5 \cr 1&2&4&5 } }}
+\tiny{{\tableau*[scY]{ 3 \cr 1&2&4&5&5 } }} \, \right) \\
& = \quad \tiny{{\tableau*[scY]{ 5 \cr 4 & 6 \cr 1&2&3 } }}
+ \tiny{{\tableau*[scY]{ 6 \cr 4 \cr 1 & 2&3&5 } }}
+\tiny{{\tableau*[scY]{ 3&4&6  \cr 1&2&5 } }}
+\tiny{{\tableau*[scY]{ 4&6 \cr 1 & 2 & 3 & 5  } }}
+\tiny{{\tableau*[scY]{ 4 \cr 1 & 2&3&5 &6  } }}  \, .
\end{split}
\end{equation}
The main goal of this section is to show that sequences of
these operators 
can be identified under $\digamma$ with the corresponding sequences of 
the operators $B_2^{(0)}$ and $B_2^{(1)}$ 
introduced in the previous section.  For that purpose, we first need to prove the
analogs of Properties 4 and  5.  
\begin{property} With $\lambda$ a partition such that $\ell(\lambda) \leq 2$, 
we have
\begin{equation}
 \mathbb B_2^{(0)} \mathbb H_{\lambda} = \mathbb H_{\lambda+1^2} \, .
\end{equation}
\end{property}
\noindent {\it Proof.} \quad 
Since $A_{n+1,n+1}$ commutes with $VS^{(n)}$, 
using the action of $A_{n+1,n+1}$ on a semi-standard tableau
of degree $n=|\lambda|$ and (4.14), we have 
\begin{equation}
A_{n+1,n+1} \sum_{\sstab ; ev(\sstab)
=\lambda'}VS^{(n)}(\sstab)
=  \sum_{\sstab ; ev(\sstab)=(\lambda',2)}VS^{(n)}(\sstab) \, ,
\end{equation}
where $\lambda'$ is a vector of length $n$.  Thus, by the
definition of $\mathbb B_2^{(0)}$,
\begin{equation}
 \mathbb B_2^{(0)} \sum_{\sstab ; ev(\sstab)=\lambda'}  VS^{(n)}(\sstab) 
 =  \sum_{\sstab ; ev(\sstab)=(\lambda',2)} 
\tau_1 r_{(11 \to 01)} \sigma_1 \cdots \sigma_{n} VS^{(n)} (\sstab)\, .
\end{equation}
The recursive nature of $VS$ gives that
$VS^{(n+2)} = \tau_1 r_{(11 \to 01)} \sigma_1 \cdots \sigma_{n} VS^{(n)} $,
which yields
\begin{equation}
\begin{split}
\mathbb B_2^{(0)} \sum_{\sstab ; ev(\sstab)=\lambda'}VS^{(n)}(\sstab) 
& = \sum_{\sstab ; ev(\sstab)=(\lambda',2)} 
VS^{(n+2)}(\sstab)  \\
& = \sum_{\bar \sstab ; ev(\bar \sstab)=(\lambda+1^2)'} 
VS^{(n+2)}(\bar \sstab) \, ,
\end{split}
\end{equation}
the last step following from  
\begin{equation}
 \sum_{\sstab ; ev(\sstab)= \mu } VS(\sstab) 
  =  \sum_{ \bar \sstab ; ev(\bar \sstab)=\beta (\mu) } VS(\bar \sstab) \, ,
\end{equation}
where $\beta(\mu)$ is any permutation of the vector $\mu$.  Formula (4.23) is a consequence
of the fact that the symmetric group action (4.6) preserves charge and shape.
Property 12 now follows from (4.22) and (4.14). 
\hfill $\square$

\begin{property}  
For $\epsilon\in\{0,1\}$, we have that 
\begin{equation}
\mathbb H_{2m+\epsilon} 
= (\mathbb B_2^{(0)} + \mathbb B_2^{(1)})^{m} \,  \mathbb H_{\epsilon} \, .
\end{equation}
\end{property}

\noindent {\it Proof. } \quad 
From Definition 10,
\begin{equation}
\mathbb  H_{n} = \sum_{\sstab ; ev(\sstab)=(1^n)} VS(\sstab) 
= \sum_{\stab} \stab \,,
\end{equation}
and
\begin{equation}
\mathbb H_{n+1,1} = 
\sum_{\sstab ; ev(\sstab)=(2,1^{n})} VS(\sstab) 
= \sum_{\stab' ; \stab'_{\{1,2\}}=12} \stab' \,
\end{equation}
where $\stab$ and $\stab'$ are standard tableaux
of degree $n$ and $n+2$ respectively.
Property 12 gives that 
$\mathbb B_{2}^{(0)}\mathbb H_n=\mathbb H_{n+1,1}$,  implying
\begin{equation}
\mathbb B_{2}^{(0)} 
\mathbb H_n =
\sum_{\stab' ; \stab'_{\{1,2\}}=12} \stab' \,.
\end{equation}
On the other hand, since $\mathbb B_2^{(1)} \stab=(\mathbb B_2^{(0)}\stab^t)^t$
and $\mathbb H_n = \mathbb H_n^t$, we have
\begin{equation}
\mathbb B_{2}^{(1)}
\mathbb H_n =
\left( \mathbb B_{2}^{(0)} \mathbb H_n \right)^t
= \sum_{\stab' ; \stab'_{\{1,2\}}=21} \stab' \,.
\end{equation}
$\mathbb H_{n+2}$ is the sum of all standard tableaux 
of order $n\!+\!2$, of which each tableaux contain
the subword 12 or 21. Therefore we have proved
\begin{equation}
\mathbb H_{n+2} = (\mathbb B_2^{(0)} + \mathbb B_2^{(1)}) \, \mathbb H_n \, ,
\end{equation}
which proves the proposition. \hfill $\square$

In order to be able to identify $\mathbb B_2^{(0)}$ and $\mathbb B_2^{(1)}$
 under $\digamma$ with $B_2^{(0)}$ and $B_2^{(1)}$, 
the main step now consists in showing that there exists
some sort of $q$-commutation relation between $\mathbb B_2^{(0)}$ and
$\mathbb B_2^{(1)}$, similar to relation (3.21) between $B_2^{(0)}$ and $B_2^{(1)}$.
This task will take most of the remainder of this section.

For any sum $S= \sum_{k} \tab^{(k)}$, where 
$\tab^{(k)} \neq \tab^{(k')}$ for all $k\neq k'$, we will say that
$\tab \in S$ if and only if 
$\tab = \tab^{(k)}$ for some $k$.  
 
We now define other linear operators on standard tableaux.
\begin{definition}
On any standard tableau $\stab$ with $\stab_{\{1,2\}}=12$ and $|\stab|=n$,
\begin{equation}
\overset{_*}{\mathbb B}_2^{(0)}  : 
\stab  \to R_{n-1} 
\sigma_{n-2}\cdots \sigma_1 r_{(01 \to 11)} \tau_{-1} \stab \,,
\end{equation}
and on any standard tableau $\stab$ such that 
$\stab_{\{1,2 \}}=21$ and $|\stab|=n$, 
\begin{equation}
\overset{_*}{\mathbb B}_2^{(1)} :
\stab  \to ( \overset{_*}{\mathbb B}_2^{(0)} \stab^t)^t \, .
\end{equation}
\end{definition}
\noindent{\it Example:} \quad Acting with $\overset{_*}{\mathbb B}_2^{(1)}$ on $\tiny{\tableau*[scY]
{ 5 \cr 2& 6 \cr 1&3&4 }}$, we go through the following steps:
\begin{equation}  
\tiny{\tableau*[scY]{ 5 \cr 2& 6 \cr 1&3&4 }} \, \overset{t}{\to} \, 
\tiny{\tableau*[scY]{ 4 \cr 3& 6 \cr 1&2&5 }} \, \overset{(1)}{\to} \, 
\tiny{\tableau*[scY]{ 3 \cr 2& 5 \cr 1&1&4 }} \, \overset{(2)}{\to} \, 
\tiny{\tableau*[scY]{ 3 \cr 2& 5 \cr 1&2&4 }} \, \overset{(3)}{\to} \, 
\tiny{\tableau*[scY]{ 3 \cr 2& 5 \cr 1&3&4 }} \, \overset{(4)}{\to} \, 
\tiny{\tableau*[scY]{ 3 \cr 2& 5 \cr 1&4&4 }} \, \overset{(5)}{\to} \, 
\tiny{\tableau*[scY]{ 3 \cr 2& 5 \cr 1&4&5 }} \, \overset{(6)}{\to} \, 
\tiny{\tableau*[scY]{ 3 \cr 2 \cr 1&4 }} \, \overset{t}{\to} \, 
\tiny{\tableau*[scY]{ 4 \cr 1&2&3 }} \, .
\end{equation}

  These operators, due to the following property,
 may be seen as inverses
of $\mathbb B_2^{(0)}$ and $\mathbb B_2^{(1)}$.
\begin{property}  
For $\stab_1'\in B_2^{(0)} \stab_1$ 
and $\stab_2'\in B_2^{(1)} \stab_2$, we have
\begin{equation}
\overset{_*}{\mathbb B}_2^{(0)}  \stab_1' = \stab_1 \quad
\hbox{and}\quad
\overset{_*}{\mathbb B}_2^{(1)} \stab_2' = \stab_2 \, . 
\end{equation}
\end{property}
\noindent {\it Proof.} \quad 
The definition of each elementary operation contained in the operators
and the fact that, for any tableau $\tab' \in A_{a,a} \tab$,
we have $R_{a} \tab'= \tab$
lead directly to the property. \hfill $\square$

\begin{property}  
The action of $\overset{_*}{\mathbb B}_2^{(1)}$
on any standard tableau $\stab$ where
$\stab_{\{ 1,2\}}=21$ and $|\stab|=n$, can equivalently be expressed as
\begin{equation}
 \overset{_*}{\mathbb B}_2^{(1)}  : 
\stab \to R_{n-1} 
\bar \sigma_{n-2}\cdots \bar \sigma_1 r_{(10 \to 11)} \tau_{-1} \stab
\end{equation} 
where $\bar \sigma_a$ is a permutation defined by its action 
on $a$ and $b=a+1$;
\begin{equation}
 a a b \overset{\bar\sigma_a}{\leftrightarrow} bab ,  
\quad  a b a  \overset{\bar\sigma_a}{\leftrightarrow} abb, \quad 
b aa  \overset{\bar\sigma_a}{\leftrightarrow} bba \, .
\end{equation}
This definition of $\bar \sigma_a$ is such that on a word $w$ with 
$ev(w_{\{a,b\}}) \in \{(1,2),(2,1) \}$, we have $\sigma_a w^R=(\bar \sigma_a w)^R$, where
the superscript $R$ stands for
the operation that 
sends a word $w=w_1 \cdots w_n$ to the word $w^R=w_n \cdots w_1$.
\end{property}
\noindent{\it Example:} \quad We have that the action (4.34) of $\overset{_*}{\mathbb B}_2^{(1)}$ on $\tiny{\tableau*[scY]{ 5 \cr 2& 6 \cr 1&3&4 }}$ goes through the following steps:
\begin{equation}  
\tiny{\tableau*[scY]{ 5 \cr 2& 6 \cr 1&3&4 }} \, \overset{(1)}{\to} \, 
\tiny{\tableau*[scY]{ 4 \cr 1& 5 \cr 1&2&3 }} \, \overset{(2)}{\to} \, 
\tiny{\tableau*[scY]{ 4 \cr 2& 5 \cr 1&2&3 }} \, \overset{(3)}{\to} \, 
\tiny{\tableau*[scY]{ 4 \cr 3& 5 \cr 1&2&3 }} \, \overset{(4)}{\to} \, 
\tiny{\tableau*[scY]{ 4 \cr 4& 5 \cr 1&2&3 }} \, \overset{(5)}{\to} \, 
\tiny{\tableau*[scY]{ 5 \cr 4& 5 \cr 1&2&3 }} \, \overset{(6)}{\to} \, 
\tiny{\tableau*[scY]{ 4 \cr 1&2&3 }} \, .
\end{equation}
Comparing with formula (4.32), we see that the transposition steps are 
avoided in (4.36) and that the tableaux obtained in (4.36)
after each steps (1 to 6) are the transposed of the corresponding tableaux in (4.32).

\noindent {\it Proof.} \quad   
For any standard tableau $\stab$ such that $\stab_{\{1,2\}}=21$ and $|\stab|=n$, we must show,
from (4.30),(4.31) and (4.34), that
\begin{equation}
\left( \sigma_{n-2}\cdots  
\sigma_1 r_{(01 \to 11)} \tau_{-1} \stab^t \right)^t 
=  \bar \sigma_{n-2}\cdots \bar \sigma_1 r_{(10 \to 11)} 
\tau_{-1} \stab \, ,
\end{equation}
where we have used the fact that $R_{n-1}$ commutes with transposition.
If we let $\bar \sstab$ be the semi-standard tableau 
$\bar \sstab =r_{(01 \to 11)} \tau_{-1} \stab^t$, we easily get that 
$\bar \sstab^t=r_{(10 \to 11)} \tau_{-1} \stab$.
Now, using again and again the identity
\begin{equation}
(\sigma_a\sstab)^t=\bar\sigma_a\sstab^t
\,\,\hbox{for any semi-standard}\,\, \sstab\,\,\hbox{with}\,\, 
ev(\sstab_{\{a,b=a+1\}})\in\{(1,2),(2,1)\}\,,
\end{equation}
we have $(\sigma_{n-2}\cdots \sigma_1 \bar \sstab)^{t}= \bar \sigma_{n-2} (\sigma_{n-3}\cdots \sigma_1 \bar \sstab)^{t}= \cdots=
\bar \sigma_{n-2} \cdots \bar \sigma_1 \bar \sstab^t$,  proving (4.37).
Identity (4.38) can be verified by observing that
under such conditions, we have
$(\sstab_{\{a,b\}})^R={\sstab_{\{a,b\}}^{t}} $.
For example, $\sstab=\tiny{\tableau*[scY]{ 4 \cr 1 &2&3&3}}$
is such that  $\sstab_{\{3,4\}}=433$ and $\sstab_{\{3,4 \}}^{t}=334$.
The only possible planar distribution of $a$ and $b$ that could
cause this to fail are
$\tiny{\tableau*[scF]{ \bl & a \cr b & \bl }}$ or 
$\tiny{\tableau*[scF]{ \bl & b \cr a & \bl }}$, since
$\sstab_{\{a,b\}}$ would be the same as $\sstab^t_{\{a,b\}}$.
The first case never holds and the second occurs only if
we have both another $a$ and another $b$, which we do not. 
Thus, using $\sstab_{\{a,b\}}^R= (\sstab^R)_{\{a,b \}}=(\sstab_{\{a,b\}})^R= 
\sstab^t_{\{a,b \}}$, we have $\bar \sigma_a \sstab_{\{ a,b\}}^t= \bar \sigma_a \sstab_{\{ a,b\}}^R= ( \sigma_a \sstab_{\{a,b \}})^R=
(\sigma_a \sstab)_{\{a,b \}}^R=(\sigma_a \sstab)_{\{ a,b\}}^t$, where the second equality 
follows from 
the definition of $\bar \sigma_a$.  This gives that
$(\sigma_{a} \sstab)^t = \bar \sigma_{a} \sstab^{t} $
on a semi-standard tableau.
\hfill $\square$

\smallskip

Let us now define an operation $\Sigma_i$ 
to act on pairs of words as follows;
\begin{equation}
\Sigma_i : (w_1,w_2) \to ( \bar \sigma_i 
\sigma_{i+1} w_1, \sigma_i \bar \sigma_{i+1} w_2)\,,
\end{equation}
and consider 6 pairs of words with $a$ and $b$ consecutive numbers;
\begin{equation}
\begin{split}
C_1(a) &= (aabb,baab), \qquad C_2(a) = (abab,abab), \qquad C_3(a) = (abba,aabb)
\\
C_4(a) &= (baab,bbaa), \qquad C_5(a) 
= (baba,baba), \qquad C_6(a) = (bbaa,abba)\,.\end{split}
\end{equation}
We insert letters into such pairs where insertion is defined with 
the operator on words,
\begin{equation}
I_k^{(a)} w_1\cdots w_n = w_1\cdots w_{k-1} a w_k\cdots w_n
\end{equation}
with the understanding that $I_{1}^{(a)} w = aw$ and $I_{n+1}^{(a)}w=wa$. 
$I_k^{(a)}$ acts on pairs of words
by $I_k^{(a)} (w_1,w_2)=(I_k^{(a)} w_1,I_k^{(a)} w_2)$.

Computer experimentation using {\em ACE}
revealed that applying $\Sigma_i$ to any pair
$C_j(i)$ with the letter $(i+2)$ inserted
recovered a pair of the same type where $i \to i+1$
with the extra letter $i$ occuring
in the same position of both elements of the pair.
More exactly,
\begin{lemma}  
Let $1\leq k,k'\leq 5$, 
$1\leq j,j'\leq 6$ and $i > 1$.
For any $k,j$ and $i$, we have that
\begin{equation}
\Sigma_i \Bigl( I_{k}^{(i+2)} C_j(i) \Bigr) = I_{k'}^{(i)} C_{j'}(i+1),
\end{equation}
for some $k'$ and $j'$.  
\end{lemma}
\noindent{\it Example:} \quad  We have
\begin{equation}
\begin{split}
\Sigma_{4} \left( I_2^{(6)} C_3(4) \right) 
& = \Sigma_4 (46554,46455) \\
& = (\bar \sigma_4 \sigma_5 46554, \sigma_4 \bar \sigma_5 46455) \\
& = (46565,46565) =  I_1^{(4)} C_5 (5) \, .
\end{split}
\end{equation}
\noindent {\it Proof.} \quad This is a property having  
30 possible configurations which are easily verified
with a computer. \hfill $\square$

\begin{proposition}  
Let $\stab_1$ be any standard tableau of degree $n$ such that ${\stab_1}_{\{ 1,2,3,4\}}$ 
is $4312$ or $3124$ 
and let $\stab_2=\stab_1^{2 \leftrightarrow 3}$, 
i.e. $\stab_2$ is obtained by permuting $2$ and $3$ in $\stab_1$.
If 
\begin{equation}
\mathbb B : \left(\stab_1,\stab_2\right) \to \Bigl(\overset{_*}{\mathbb B}_2^{(1)}
\overset{_*}{\mathbb B}_2^{(0)}\stab_1 ,
\overset{_*}{\mathbb B}_2^{(0)} 
\overset{_*}{\mathbb B}_2^{(1)} \stab_2 \Bigr)\,,
\end{equation}
then 
\begin{equation}
\mathbb B \left(\stab_1,\stab_2\right) = (\bar \stab, \bar \stab)\, ,
\end{equation}
for some standard tableau $\bar \stab$.
\end{proposition}

\noindent{\it Example:} \quad 
Given $\stab_1=\tiny{\tableau*[scY]{ 7 \cr 3&5 \cr 1&2&4&6&8 }}$, 
we have the pair, 
\begin{equation}
\mathbb B
(\stab_1,\stab_2)=
\mathbb B \left(\,\tiny{\tableau*[scY]{ 7 \cr 3&5 \cr 1&2&4&6&8 }}\,,
\tiny{\tableau*[scY]{ 7 \cr 2&5 \cr 1&3&4&6&8 }} \, \right)\,
=\left( \, \tiny{\tableau*[scY]{  2 \cr 1&3&4 }} \, , \, 
\tiny{\tableau*[scY]{  2 \cr 1&3&4 }} \, \right)\, .
\end{equation}

\noindent {\it Proof.} \quad
$\overset{_*}{\mathbb B}_2^{(0)}\!$ reduces the degree of 
$\stab_1$ to $n\!-\!2$, giving by Definition 14 and Property 16,
\begin{equation}
\begin{split}
\overset{_*}{\mathbb B}_2^{(1)} \overset{_*}{\mathbb B}_2^{(0)}\! \stab_1 & = 
R_{n-3} \bar \sigma_{n-4} \cdots \bar \sigma_1
r_{(10 \to 11)} \tau_{-1} 
R_{n-1}   \sigma_{n-2} \cdots  \sigma_1 r_{(01 \to 11)} \tau_{-1}
\stab_1 \\
 =&R_{n-3}R_{n-2} \bigl( \bar \sigma_{n-4}  \sigma_{n-3} \bigr) \cdots 
\bigl(\bar \sigma_{1} \sigma_2 \bigr)  
r_{(10 \to 11)} \tau_{-1}  \sigma_2  \sigma_1 r_{(01 \to 11)} \tau_{-1} \stab_1\,, 
\end{split}
\end{equation}
where we have considered the relations $\tau_{-1}R_{n-1}=R_{n-2}\tau_{-1}$ and 
$\tau_{-1}\sigma_i=\sigma_{i-1}\tau_{-1}$.
Similarly, acting first with $\overset{_*}{\mathbb B}_2^{(1)}\!$,
\begin{equation}
\begin{split}
\overset{_*}{\mathbb B}_2^{(0)} \overset{_*}{\mathbb B}_2^{(1)}\! \stab_2
&= R_{n-3}  \sigma_{n-4} \cdots  \sigma_1
r_{(01 \to 11)} \tau_{-1} R_{n-1} \bar  \sigma_{n-2} \cdots  \bar \sigma_1
r_{(10 \to 11)} \tau_{-1}\stab_2 \\
=&R_{n-3}R_{n-2} \bigl( \sigma_{n-4} \bar \sigma_{n-3} \bigr) 
\cdots \bigl( \sigma_{1}
\bar \sigma_2 \bigr)  
r_{(01 \to 11)} \tau_{-1} \bar \sigma_2 \bar \sigma_1 
r_{(10 \to 11)} \tau_{-1} \stab_2 \,.
\end{split}
\end{equation}
We act first with $r_{(10 \to 11)}\tau_{-1}\sigma_2\sigma_1r_{(01 \to 11)}\tau_{-1}$ 
on $\stab_1 $ in the case that ${\stab_1}_{\{1,2,3,4 \}}$ is $3124$,
obtaining a word with subword $1122$ and letters $3,\dots,n\!-\!2$ 
occurring exactly once. 
$\stab_2$, defined by permuting $2$ and $3$ in $\stab_1$, 
thus contains the subword $2134$ which is sent to $2112$ under 
$r_{(01 \to 11)}\tau_{-1}\bar\sigma_2\bar\sigma_1 r_{(10\to 11)}\tau_{-1}$,
while the remaining letters occur exactly as they do in 
$r_{(10 \to 11)}\tau_{-1}\sigma_2\sigma_1r_{(01 \to 11)}\tau_{-1} \stab_1$.
Consequently,
\begin{equation}
\mathbb B \left( \stab_1, \stab_2 \right) = 
R_{n-3}R_{n-2} \Sigma_{n-4} \cdots \Sigma_{1} 
\left(I_{k_1}^{(n-2)} \cdots I_{k_{n-5}}^{(4)}I_{k_{n-4}}^{(3)} 
(1122,2112)\right)\,,
\end{equation}
for some $k_1,\dots,k_{n-4}$.
Since $\Sigma_1$ acts only on the letters 1,2 and 3, we may now use Lemma 17, where $i=j=1$, 
to determine that the 
action of $\Sigma_1$ results in 
\begin{equation}
\mathbb B \left( \stab_1, \stab_2 \right) = 
R_{n-3}R_{n-2} \Sigma_{n-4} \cdots \Sigma_{2} \left(I_{k_1}^{(n-2)}
\cdots I_{k_{n-5}}^{(4)}I_{r}^{(1)} C_j(2)\right)\,,
\end{equation}
for some $j$ and $r$.  
Lemma 17, applied repeatedly in this manner, gives
\begin{equation}
\mathbb B \left( \stab_1, \stab_2\right) = 
R_{n-3}R_{n-2}\left(I_{r_1}^{(1)}I_{r_2}^{(2)}\cdots I_{r_{n-4}}^{(n-4)}C_j'(n-3)\right)\,,
\end{equation}
for some $r_1,\dots,r_{n-4}$ and some $j'$.  
This is to say that the action of $\mathbb B$ on such a pair 
is equivalent to acting with 
$R_{n-3}R_{n-2}$ on a pair of tableaux that are
identical in all letters except $n\!-\!3$ and $n\!-\!2$. 
Further, since $R_{n-3}R_{n-2}$ removes these letters,
we have proved the identity in the case $3124$.  
A sequence of similar arguments may be used in the case 
that $ {\stab_1}_{\{1,2,3,4 \}}$ is $4312$, and we get 
\begin{equation}
\mathbb B \left( \stab_1, \stab_2\right) = 
R_{n-3}R_{n-2} \Sigma_{n-4} \cdots \Sigma_{1} 
\left(I_{k_1}^{(n-2)} \cdots I_{k_{n-5}}^{(4)}I_{k_{n-4}}^{(3)} 
(2112,2211)\right)\,.
\end{equation}
Again, successive applications of Lemma 17, beginning with
the case $i=1$ and $j=4$, prove the identity.
\hfill$\square$

\begin{lemma} 
For $\stab'$ a standard tableau of degree $n \geq 4$, we have that
\begin{equation}
\begin{split}
\stab'\in \mathbb B_2^{(0)} \mathbb B_2^{(0)}\sum_{\stab; |\stab|=n-4} \stab 
&\iff \stab'_{\{1,2,3,4\}}\in\{1234, 4123,3412\}\,, \\
\stab'\in \mathbb B_2^{(0)} \mathbb B_2^{(1)}\sum_{\stab ; |\stab|=n-4} \stab 
&\iff \stab'_{\{1,2,3,4 \}}\in\{4312,3124\}\,,\\
\stab'\in \mathbb B_2^{(1)} \mathbb B_2^{(0)}\sum_{\stab ; |\stab|=n-4} \stab 
&\iff \stab'_{\{1,2,3,4 \}}\in\{4213,2134\}\,,\\
\stab'\in \mathbb B_2^{(1)} \mathbb B_2^{(1)}\sum_{\stab ; |\stab|=n-4} \stab 
&\iff \stab'_{\{ 1,2,3,4\}}\in\{4321, 3214,2413\}\,.\\
\end{split}
\end{equation}
\end{lemma}
\noindent {\it Proof.} 
\quad 
We begin by simultaneously proving the first
two cases of ($\Rightarrow$) and the 
others of ($\Rightarrow$) follow by transposition.
Notice first that $B_{2}^{(0)} \stab=\sum_{\stab''} \stab''$, 
where $\stab''_{\{1,2\}}$ is $12$ and
$B_{2}^{(1)} \stab= \sum_{\stab''} \stab''$, 
where $\stab''_{\{1,2 \}}$ is $21$ (see the proof of Property 13). 
The following action of $B_{2}^{(0)}$ 
begins with $A_{n+1,n+1}$ adding a horizontal 
$2$ strip to $\stab''$ resulting in tableaux that are all 
semi-standard and containing the subword $12$ (or $21$).  
We act next with the succession of  
$\sigma_n,\sigma_{n-1},\dots,\sigma_3$
implying that the tableaux remain semi-standard
and thus must each contain, 
in the first case, the subword $1233$, $3123$ or $3312$,
and in the second,  $3213$ or $2133$. 
The remaining operations, aside from $\tau_1$,
act exclusively on these subwords as follows;
\begin{equation}
\begin{split}
\tau_1 r_{(11 \to 01)} \sigma_1 \sigma_2 1233 &= 1234\,,\\
\tau_1 r_{(11 \to 01)} \sigma_1 \sigma_2 3123 &= 4123\,,\\
\tau_1 r_{(11 \to 01)} \sigma_1 \sigma_2 3312 &= 3412\,,\\
\end{split}
\end{equation}
for the first case and 
\begin{equation}
\begin{split}
\tau_1 r_{(11 \to 01)} \sigma_1 \sigma_2 3213 &= 4312\,,\\
\tau_1 r_{(11 \to 01)} \sigma_1 \sigma_2 2133 &= 3124\,,\\
\end{split}
\end{equation}
for the second, thus proving $(\Rightarrow)$ for the two first cases.  
To prove $(\Leftarrow)$, we are given a standard tableau 
$\stab'$ with subword $\stab'_{\{1,2,3,4\}}$ in one of the
four defined disjoint sets; call this set $S_{\epsilon_1,\epsilon_2}$. 
Property~13 gives that $\stab'$, which is $\in \mathbb H_n$, for  
$n=|\stab'|$,  
is such that $\stab'\in \mathbb B_2^{(\bar\epsilon_1)}\mathbb B_2^{(\bar\epsilon_2)}
\mathbb H_{n-4} =\mathbb B_2^{(\bar\epsilon_1)}\mathbb B_2^{(\bar\epsilon_2)}
\sum_{\{|\stab|=n-4\}}\stab $, for some $\bar \epsilon_i \in \{0,1 \}$.
But since we have just proven that for such
$\stab'$,
$\stab'_{\{1,2,3,4\}}$ is contained in the
set $S_{\bar\epsilon_1,\bar\epsilon_2}$,
we see that $\bar\epsilon_i=\epsilon_i$ and 
the proposition is proven.  \hfill$\square$
\smallskip

\begin{proposition} 
On any standard tableau $\stab$, we have
\begin{equation}
\mathbb B_2^{(0)} \mathbb B_2^{(1)} \stab= 
\left( \mathbb B_2^{(1)} \mathbb B_2^{(0)} 
\stab \right)^{2 \leftrightarrow 3} \, ,
\end{equation}
where $2\!\leftrightarrow\!3$ denotes a permutation of the 
letters $2$ and $3$ in each tableaux.
\end{proposition}
\noindent {\it Proof.}  \quad 
Suppose there exists $\stab'\!\in\!\mathbb B_2^{(0)}\mathbb B_2^{(1)}\stab$ 
such that 
${\stab'}^{2\leftrightarrow 3}\!\not\in\!
\mathbb B_2^{(1)}\mathbb B_2^{(0)}\stab$. 
Lemma 19 gives that every element in 
$\mathbb B_2^{(0)}\mathbb B_2^{(1)} \stab$,
in particular $\stab'$, 
must contain either the subword $4312$ or $3124$.  This
implies that ${\stab'}^{2\leftrightarrow 3}$ must contain either 
$4213$ or $2134$ and thus by the same lemma we have that 
\begin{equation}
{\stab'}^{2\leftrightarrow 3}\in \mathbb B_2^{(1)}\mathbb B_2^{(0)} \stab''
\end{equation}
for some $ \stab''\neq \stab$. 
Observe now that
$\overset{_*}{\mathbb B}_2^{(1)}\overset{_*}{\mathbb B}_2^{(0)}\stab'= \stab$, 
implies by Proposition 18 that 
$\overset{_*}{\mathbb B}_2^{(0)}\overset{_*}{\mathbb B}_2^{(1)}
{\stab'}^{2\leftrightarrow 3}$
must also be $\stab$.  
Expression (4.57) then yields
\begin{equation}
\overset{_*}{\mathbb B}_2^{(0)}\overset{_*}{\mathbb B}_2^{(1)}
\left(\mathbb B_2^{(1)}\mathbb B_2^{(0)} \stab'' \right)
=\stab + \hbox{other terms}\,
\end{equation}
which by Property 15 gives that $\stab''=\stab$ 
and we reach a contradiction.  We thus have that $\mathbb B_2^{(0)} \mathbb B_2^{(1)} \stab 
\subseteq  
\left( \mathbb B_2^{(1)} \mathbb B_2^{(0)} 
\stab \right)^{2 \leftrightarrow 3}$.  We can also show in the same manner that
$\mathbb B_2^{(0)} \mathbb B_2^{(1)} \stab 
\supseteq  
\left( \mathbb B_2^{(1)} \mathbb B_2^{(0)} 
\stab \right)^{2 \leftrightarrow 3}$, which proves the proposition.
\hfill$\square$
\smallskip

We now define four pairs of words on
consecutive numbers, $a,b,c$ and $d$;
\begin{equation}
\begin{split}
& D_1(a) = ( b a c d, c a b d), \qquad \qquad D_2(a) = (dbac, dcab),  \\ 
& D_3(a) = (a c d b, abdc), \qquad \qquad
D_4(a) = (cd ba, bdca)\, .
\end{split}
\end{equation}
These pairs often appear as the only distinct
subwords in a given pair of semi-standard tableaux.  More
precisely, such a pair of semi-standard tableaux, 
called $(\sstab_1,\sstab_2)_{D_j(a)}\,$, satisfies
$\sstab_1=\sstab_2^{b\leftrightarrow c}$ 
and $(\sstab_{1_{\{a,b,c,d\}}},\sstab_{2_{\{a,b,c,d\}}})=D_j(a)$.
For example, 
\begin{equation}
\left( \, \tiny{\tableau*[scY]{  5 \cr 4&7&8 \cr 1&2&3&6 }} \, , \, 
\tiny{\tableau*[scY]{  5 \cr 4&6&8 \cr 1&2&3&7 }} \, \right)_{D_3(5)} 
\end{equation}
is such a pair.  One should note that $(\sstab_1,\sstab_2)_{D_j(a)}$ 
is a pair of tableaux of the same shape since 
in any such semi-standard tableaux, $b$ and $c$ 
never occur in the same row or column.
With $\Omega_i$ defined such that on pairs of words
\begin{equation}
\Omega_i : (w_1,w_2) \to (\sigma_i \sigma_{i+1} \sigma_{i+2} 
\sigma_{i+3} w_1, \sigma_i \sigma_{i+1} \sigma_{i+2} \sigma_{i+3} w_2)\,;
\end{equation}

\begin{lemma} Let $1 \leq k_2,k_2',k_2'' \leq 5$, \,
$1 \leq k_1,k_1', k_1'' \leq 6$,\,  $1 \leq j,j' , j''\leq 4$ and $i>0$.
For any such $k_1,k_2,j$ and $i$, we have
\begin{equation} 
r_{(ii \to i i-1)} \Omega_i I_{k_1}^{(i+4)} I_{k_2}^{(i+4)} D_{j}(i)
\in\left\{ 
I_{k_1'}^{(i-1)} I_{k_2'}^{(i)} D_{j'}(i+1)\,, \,
I_{k_1'}^{(i-1)}I_{k_2'}^{(i+4)} D_{j'}(i) \right\}
\end{equation}
and
\begin{equation} 
r_{(ii \to i-1 i)} \Omega_i I_{k_1}^{(i+4)} I_{k_2}^{(i+4)} D_{j}(i) 
\in\left\{I_{k_1''}^{(i-1)} I_{k_2''}^{(i)} D_{j''}(i+1)\,, \,
I_{k_1''}^{(i-1)}I_{k_2''}^{(i+4)} D_{j''}(i)\right\}\, ,
\end{equation}
for some $k_1',k_1'', k_2',k_2'', j'$ and $j''$. 
\end{lemma}
\noindent{\it Example:} \quad 
Starting with $I_4^{(9)}I_2^{(9)}D_4(5)=(798965,698975)$, we get
\begin{equation}
\Omega_4 (798965,698975) = (597865,596875) \,.
\end{equation}
Under $r_{(55 \to 45)}$, we recover
$I_1^{(4)} I_1^{(9)} D_4(5)$ and under
$r_{(55 \to 54)}$, $I_6^{(4)} I_2^{(9)} D_3(5)$.

\noindent {\it Proof.} \quad For each $i$ there are 60 cases that 
have been verified using a computer.\hfill$\square$

\begin{lemma} 
Let $\stab$ and $\stab'$ be standard tableaux
of type $(\stab,\stab')_{D_j(i)}$ for some $i,j$. 
If standard tableaux $\bar \stab \in \mathbb B_2^{(\epsilon)} \stab$ and
$\bar \stab'\in \mathbb B_2^{(\epsilon)} \stab'$
are standard tableaux of the same shape for $\epsilon\in\{0,1\}$,  
then $\bar \stab$ and $\bar \stab'$ is a pair of 
type $(\bar \stab,\bar \stab')_{D_{j'}(i+1)}$ or 
$(\bar \stab,\bar \stab')_{D_{j'}(i+2)}$, for some $j'$.
\end{lemma}

\noindent {\it Proof.} \quad We start 
with the case $\epsilon = 0$ and split 
the action of $B_2^{(0)}$ into a sequence of
operations beginning with $A_{n+1,n+1}$.
As such, we consider a semi-standard tableau  $\sstab$
obtained by adding an arbitrary horizontal $2$-strip to $\stab$.
We denote by $\sstab'$, the semi-standard tableau of the same shape that is 
obtained by adding this horizontal $2$-strip to $\stab'$
and thus $\sstab$ and $\sstab'$ are a pair of type $(\sstab,\sstab')_{D_j(i)}$. 
Next in the sequence of operations defining $B_2^{(0)}$ is
$\sigma_{i+4}\cdots\sigma_n$ which, acting
only on the letters $i\!+\!4,\dots,n$, must
preserve the similarity in $\sstab$ and $\sstab'$.
If $i \neq 1$, since acting with $\sigma_{i-1}$ amounts to applying either
$r_{(ii \to i i-1)}$ or $r_{(ii \to i-1 i)}$,
acting on both elements with 
$\sigma_{i-1}\sigma_i\sigma_{i+1}\sigma_{i+2}\sigma_{i+3}$ using
Lemma 21, gives a pair of semi-standard tableaux of type 
$(\bar \sstab,\bar \sstab')_{D_{j'}(i+1)}$ or 
$(\bar \sstab,\bar \sstab')_{D_{j'} (i)}$.
There remains to act with $\tau_1r_{11\to 01}\sigma_1\dots\sigma_{i-2}$, which leads to
 pairs of standard tableaux of type 
$(\bar \stab,\bar \stab')_{D_{j'}(i+2)}$ or 
$(\bar \stab,\bar \stab')_{D_{j'} (i+1)}$.  In the case where $i=1$, 
acting on both elements with $\tau_1 r_{(11 \to 0 1)} \sigma_1\sigma_{2}\sigma_{3}\sigma_{4}$,
gives, from Lemma 21, pairs of standard tableaux of type 
$(\tilde \stab,\tilde \stab')_{D_{j'}(3)}$ or 
$(\tilde \stab,\tilde \stab')_{D_{j'} (2)}$, 
finally proving the lemma for $\epsilon=0$.  
To prove the lemma in the case 
$\mathbb B_{2}^{(1)}$, 
we observe that, for $\stab_1$ and $\stab_2$ standard tableaux,
\begin{equation}
\begin{split}
& (\stab_1,\stab_2)_{D_{1}(a)} \implies (\stab_2^t,\stab_1^t)_{D_{2}(a)} \,
\qquad
(\stab_1,\stab_2)_{D_{2}(a)} \implies (\stab_2^t,\stab_1^t)_{D_{1}(a)} \, ,\\ 
& (\stab_1,\stab_2)_{D_{3}(a)} \implies (\stab_2^t,\stab_1^t)_{D_{4}(a)} \, 
\qquad
(\stab_1,\stab_2)_{D_{4}(a)} \implies (\stab_2^t,\stab_1^t)_{D_{3}(a)} \, .\\ 
\end{split}
\end{equation} 
$\mathbb B_{2}^{(1)} \stab=(\mathbb B_2^{(0)} \stab^t)^t$ thus 
implies that the proof in this case is exactly the proof 
for $\mathbb B_2^{(0)}$ with every pair reversed 
plus an additional reversal of the pairs at the end, accounting for the
last transposition in $\mathbb B_2^{(1)}$. 
\hfill $\square$ 

\begin{lemma}  If $\stab$ and $\stab'$ are 
a pair of standard tableaux of type $(\stab,\stab')_{D_j(i)}$
then
\begin{equation}
\digamma(\stab) = q \digamma(\stab') \,.
\end{equation}
\end{lemma}
\noindent {\it Proof.} \quad 
We have already noted that $(\stab,\stab')_{D_j(i)}$ 
is a pair of tableaux with the same shape.
Further, the definition of cocharge gives
$\cocharge(\stab)=\cocharge(\stab')+1$; 
for example, for $D_1(a)=(bacd,cabd)$, we have
\begin{equation}
\begin{split}
\cocharge(\stab)&=c_1 + \dots + c_a + (c_a+1) + (c_a+1) +(c_a+1) + c_{a+4} 
+ \cdots +c_n\\
\cocharge (\stab')&= c_1 + \dots + c_a + (c_a) + (c_a+1) 
+(c_a+1) + c_{a+4} + \cdots +c_n\, ,
\end{split}
\end{equation} 
giving $\cocharge(\stab)=\cocharge(\stab')+1$ as claimed.
\hfill $\square$

\begin{definition}
We define, for $v=(v_1,\dots,v_k)$ with $v_i \in \{0,1 \}$, 
\begin{equation}
\mathbb U_v^{(\epsilon)} =  \mathbb B_2^{(v_1)} \cdots \mathbb B_2^{(v_k)}\mathbb  H_\epsilon \, , \qquad \epsilon \in \{0,1\} \, .
\end{equation}
\end{definition}

\begin{proposition} For any $v=(v_1,\dots,v_k)$ and 
$\bar v= (\bar v_1,\dots,\bar v_{k'})$,  with $v_i,\bar v_i \in \{0,1 \}$ and
$k,k' \geq 0$,
we have
\begin{equation}
\digamma( \mathbb U_{v,1,0,\bar v}^{(\epsilon)}) = 
q \digamma( \mathbb U_{v,0,1,\bar v}^{(\epsilon)}) \qquad \epsilon \in \{ 0,1\} \, .
\end{equation}
\end{proposition}
\noindent {Proof.} \quad 
We begin by showing that this identity holds in
the case that $v$ is empty. 
Since 
$\mathbb U_{1,0,\bar v}^{(\epsilon)} = 
\mathbb B_2^{(1)} \mathbb B_2^{(0)}\sum_{\stab} \stab$
and
$\mathbb U_{0,1,\bar v}^{(\epsilon)} = 
\mathbb B_2^{(0)} \mathbb B_2^{(1)}\sum_{\stab} \stab$,
each tableau $ \stab_1\! \in\! \mathbb U_{1,0,\bar v}^{(\epsilon)}$ 
can be paired with some $ \stab_1'\! \in\! \mathbb U_{0,1,\bar v}^{(\epsilon)}$ 
such that $\stab_1'\!=\!  \stab_1^{2\leftrightarrow 3}$ by Proposition 20
and such that ${\stab_1}_{\{1,2,3,4\}}\!\in\! \{4213,2134\}$ by 
Lemma 19.  This implies that $ \stab_1$ and $ \stab_1'$ are of type
$(\stab_1, \stab_1')_{D_j(1)}$ where $j=1$ or $2 $  which, using 
Lemma 23, proves the identity for $v=()$. 
We now proceed to the case when $v=(0)$ or $(1)$ 
by acting with $\mathbb B_2^{(\epsilon)}$
on the pairs obtained when $v=()$.
These pairs $(\stab_1, \stab_1')_{D_j(1)}$ 
are thus sent to a pair of sums of standard tableaux
that, by Lemma 22, can be paired
by types $(\stab_2 ,  \stab_2')_{D_j(i)}$ 
for some $j=1,2,3$ or $4$ and some $i=3$ or $4$. 
For any $v$, we repeat this process and obtain
that each $ \bar \stab \in \mathbb U_{v,1,0,\bar v}^{(\epsilon)}$ 
can be paired with some standard tableau 
$\bar \stab'\in\mathbb U_{v,0,1,\bar v}^{(\epsilon)}$ 
where this pair is of type $(\bar \stab,\bar \stab')_{D_{j}(i)}$, 
for $1\leq j\leq 4$ and $1\leq i\leq n-3$.  
Lemma 23 then proves the theorem.
\hfill $\square$

\begin{theorem}  Let $\epsilon,v_i \in \{0,1\}$. 
For any $v=(v_1,\dots,v_k)$ we have 
\begin{equation}
\digamma( \mathbb U_v^{(\epsilon)}) = U_v^{(\epsilon)} \,.
\end{equation}
\end{theorem} 
\noindent {\it Proof. } \quad  
Recall that the action of 
$B_2^{(0)}$ and $B_2^{(0)}+B_2^{(1)}$
was determined in Properties 4 and 5, respectively. 
We have further, by Corollary 9, that
\begin{equation}
H_{2m +\ell + \epsilon,\ell}[X;q,t] = \sum_{\bar v} U_{\bar v}^{(\epsilon)}\,,
\end{equation} 
where $\bar v=(0^{\ell},v)$ for some $v=(v_1,\dots,v_m)$.  
Observe that we have proved equivalent actions
for the operators $\mathbb B_2^{(0)}$ and 
$\mathbb B_2^{(0)}+\mathbb B_2^{(1)}$
in Properties 12 and 13, giving
\begin{equation}   
\mathbb H_{2m+ \ell + \epsilon,\ell} = \sum_{\bar v}  \mathbb  U_{\bar v}^{(\epsilon)} \,,
\end{equation}
where $\bar v$ is as before. 
We thus have, since $\digamma(\mathbb H_{m,\ell})=H_{m,\ell}[X;q,t]$,
\begin{equation}
H_{2m + \ell + \epsilon,\ell}[X;q,t] = 
\sum_{\bar v} \digamma(\mathbb  U_{\bar v}^{(\epsilon)})
=\sum_{\bar v}  U_{\bar v}^{(\epsilon)}\, .
\end{equation}
We convert the expression such that we are summing only 
over dominant vectors 
$v_d\!=\!(0^{m+\ell-k},1^k)$ for some $k$
by using the following implication of the '$q$-commutation' relations 
proven in Property 6 and Proposition 25:
$ U_{\beta(v)}^{(\epsilon)} = q^{\ell(\beta)}  U_{v_d}^{(\epsilon)}$ and
$\digamma(\mathbb U_{\beta(v)}^{(\epsilon)}) = q^{\ell(\beta)} 
\digamma(\mathbb U_{v_d}^{(\epsilon)})$,
where $\ell(\beta)$ is the length of the permutation $\beta$ such that $\beta(v)=v_d$.
This gives
\begin{equation}
H_{2m+ \ell +\epsilon,\ell}[X;q,t] =\sum_{v_d} d_{v_d}^{m,\ell}(q) \digamma(\mathbb U_{v_d}^{(\epsilon)})
= \sum_{v_d} d_{v_d}^{m,\ell}(q) U_{v_d}^{(\epsilon)}\, ,
\end{equation}
where $d_{v_d}^{m,\ell}(q)=\sum_{\beta( v) =v_d} q^{\ell(\beta)}$.
For $2m+2\ell+\epsilon=n$, the
number of possible $v_d$ is $\lfloor n/2 \rfloor +1$,
exactly the number of partitions of $n$ of length $ \leq 2$.  
We thus have, from (4.74), that
$U_{v_d}^{(\epsilon)}$  and $\digamma(\mathbb U_{v_d}^{(\epsilon)})$ 
are both bases for $\mathcal V$, the $\mathbb Q[q,t]$-linear span of $\{H_{\lambda}[X;q,t]\}_{\ell(\lambda) \leq 2}$.  
We see from expression (4.74) again, that the transition matrices from 
$\{H_{\lambda}[X;q,t]\}_{\ell(\lambda) \leq 2}$ 
to $\{ U_{v_d}^{(\epsilon)}\}_{v_d}$ and from 
$\{H_{\lambda}[X;q,t]\}_{\ell(\lambda) \leq 2}$ 
to $\{ \digamma( \mathbb U_{v_d}^{(\epsilon)})\}_{v_d}$ are identical.  
Since these are invertible matrices, we have that
$U_{v_d}^{(\epsilon)}= \digamma(\mathbb U_{v_d}^{(\epsilon)})$, which can be
extended to $U_{v}^{(\epsilon)}=\digamma(\mathbb U_{v}^{(\epsilon)})$ 
using Property 6 and Proposition 25. \hfill $\square$

\section{A statistic for Macdonald polynomials in 2 parts}

It is now clear from the previous theorem and (3.28) that for $d=2m+2\ell+\epsilon$,  we have
\begin{equation}
J_{2m+ \ell+\epsilon,\ell}[X;q,t] =  \sum_{v= (v_1,\dots,v_{m+\ell}) } 
q^{(1-d)  |v|_{\ell}+2 n(v)_{\ell}} t^{\ell -|v|_{\ell} } \digamma ( \mathbb U_v^{(\epsilon)}),
\end{equation}
where $|v|_{\ell}$ and $n(v)_{\ell}$ are as defined in Proposition~8.  
To provide an expression for $J_{\lambda}[X;q,t]$ 
with coefficients that are determined by statistics, 
we associate to any standard tableau a vector $\in \{0,1\}^k$
called a "domino" vector. 
The domino vector is determined by the succession of applications 
of $\mathbb B_2^{(0)}$ and $\mathbb B_2^{(1)}$ 
that build the associated standard tableau.
Since any standard tableau $\stab$ such that 
$\stab_{\{1,2 \}}=12$ (or 21) can be obtained by acting 
with $\mathbb B_2^{(0)}$ (or $\mathbb B_2^{(1)}$)
on a predecessor $\stab'$, 
such a succession is determined recursively
using Property 15.
\begin{theorem} The Macdonald polynomials indexed by 
partitions with $\leq 2$ parts are 
\begin{equation}
J_{2m+\ell+\epsilon ,\ell} = 
\sum_{|\stab|=d}\stat (\stab) S_{\shape(\stab)}[X^t],
\end{equation}
where $d=2m + 2 \ell+ \epsilon$ and  
\begin{equation}
\stat (\stab) = q^{\cocharge(\stab)} q^{(1-d)  |\dv(\stab)|_{\ell}+2 n(\dv(\stab))_{\ell}} t^{\ell -|\dv(\stab)|_{\ell} },
\end{equation}
with the domino vector,
$\dv(\stab)=(\dv_1,\dots,\dv_{m+\ell})$, obtained recursively by 
\begin{equation}
\dv(\stab) = 
\begin{cases}
\Bigl(0,\dv \bigl( \overset{_*}{\mathbb B}_2^{(0)} \stab \bigr) \Bigr) 
& {\text{if $\stab_{\{1,2\}}=12$} }  \\
\Bigl(1,\dv \bigl( \overset{_*}{\mathbb B}_2^{(1)} \stab \bigr) \Bigr) & 
{\text {if $\stab_{\{1,2\}}=21$ }}\\
\qquad\emptyset  & {\text {if $\stab$ has degree $\leq 1$ }}\\
\end{cases} \, .
\end{equation}
\end{theorem}

\noindent {\it Proof.} \quad  The theorem follows directly from 
(5.1) and Property 15. \hfill $\square$

\noindent{\it Example:} 
The statistic associated to a standard tableau  
$\stab=  \tiny{{\tableau*[scY]{ 4&8 \cr 3 & 5 & 7  \cr 1&2&6 } }}$
in the Macdonald polynomial $J_{6,2}[X;q,t]$ is 
determined by finding the domino vector of $\stab$.
\begin{equation}
\begin{split}
(0,\dv \left( \,  \overset{_*}{\mathbb B}_{2}^{(0)}  \, 
\tiny{{\tableau*[scY]{ 4&8 \cr 3 & 5 & 7  \cr 1&2&6 } }} 
= \tiny{{\tableau*[scY]{ 4\cr 2 & 6 \cr 1 & 3 & 5 } }} \, \right) \, ) & =
 (0,1,\dv \left( \overset{_*}{\mathbb B}_{2}^{(1)} \, 
\tiny{{\tableau*[scY]{ 4\cr 2 & 6 \cr 1 & 3 & 5 } }} 
= \tiny{{\tableau*[scY]{ 4\cr 3  \cr 1 & 2 } }}  \, \right) \, )\\
&=  (0,1, 0,\dv \left( \overset{_*}{\mathbb B}_{2}^{(0)} \, \tiny{{\tableau*[scY]{ 4\cr 3 \cr 1  & 2 } }} = \tiny{{\tableau*[scY]{  2 \cr 1  } }} \, \right) \, ) \\
& = (0,1,0,1) \, .
\end{split}
\end{equation}
This gives that $|\dv(\stab)|_2=1$ and $n(\dv(\stab))_2=1$.  
The cocharge of $\stab=48~357~126$ is $0+0+1+2+2+2+3+4=14$, and we have
\begin{equation}
\stat (\stab) = q^{14} q^{-7+2} t^{2-1}=q^9 t \, .
\end{equation}

\begin{acknow}
The authors would like to thank Alain Lascoux for his
significant contribution to this work
and for unceasingly devoting himself to answering a constant 
barrage of questions. Undoubtedly this work would
not exist without him.  This work was carried out 
at Marne-la-Vall\'ee and the authors are grateful for
the resources provided there.  The author, J. Morse, expresses thanks to Adriano Garsia for 
the NSF support that made her visit to Marne-la-Vall\'ee possible. 
L. Lapointe is supported through an NSERC postdoctoral fellowship (Canada).
The algebraic combinatorics environment, ACE,
Maple library \cite{[V]} has been used extensively.  Finally,  we thank Glenn Tesler for
allowing us to use his  tableau \LaTeX \, package. 
\end{acknow}


\begin{thebibliography}{33}
\bibitem{[Fish]} S. Fishel, \emph{Statistics for special Kostka polynomials}, Proc. Amer. Math. Soc. {\bf 123} (1995), 2961--2969.
\bibitem{[GH]} A.M. Garsia and M. Haiman, \emph{A graded representation module for Macdonald's
polynomials}, Proc. Natl. Acad. Sci. USA V {\bf 90} (1993) 3607--3610.
\bibitem{[KN]} A. Kirillov and M. Noumi,  \emph{Affine Hecke algebras and raising operators for Macdonald
polynomials}, Duke Math. J., to appear.
\bibitem{[LLM]} L. Lapointe, A. Lascoux and J. Morse, 
\emph{Determinantal expressions for Macdonald polynomials}, 
IMRN {\bf 18} (1998), 957--978.
\bibitem{[LV]}   
L. Lapointe and L. Vinet, \emph{A short proof of the integrality of the Macdonald $(q,t)$-Kostka coefficients},
Duke Math. J. {\bf 91} (1998), 205--214.
\bibitem{[LS1]} A. Lascoux and M.-P. Sch\"utzenberger, \emph{ Le mono\" \i de
plaxique}, Quaderni della Ricerca scientifica {\bf 109} (1981), 129--156.
\bibitem{[LS2]} A. Lascoux and  M.-P. Sch\"utzenberger, \emph{Sur une
conjecture de H.O. Foulkes}, C.R. Acad. Sc. Paris. {\bf 294} (1978), 323--324.
\bibitem{[Ma]} I.~G. Macdonald, \emph{Symmetric functions and Hall
    polynomials}, 2nd edition, Clarendon Press, Oxford, 1995.
\bibitem{[V]} {S. Veigneau}, {\it ACE, an Algebraic
                      Combinatorics Environment for the computer
                      algebra system MAPLE\/}, {\it Version 3.0\/},
                      Universit\'e de Marne-la-Vall\'ee, 1998,
http://phalanstere.
univ-mlv.fr/{\~{}}ace/.  
\end{thebibliography}
\end{document}